\documentclass[a4paper, 11pt]{amsart}
\usepackage{amssymb}
\usepackage{amsthm}
\usepackage{fullpage}
\usepackage{graphicx}
\usepackage{subfig}
\usepackage{amsmath,amscd}
\usepackage{xypic}
\usepackage{color}
\usepackage{float}
\usepackage{vmargin}
\usepackage[colorlinks,citecolor=blue,filecolor=black,linkcolor=blue,urlcolor=black]{hyperref}
\usepackage[all,cmtip]{xy}

\newcommand{\Z}{\ensuremath{\mathbb Z}}

\newcommand{\R}{\ensuremath{\mathbb R}}

\newcommand{\CP}{\ensuremath{\mathbb {CP}}}

\theoremstyle{plain}
\newtheorem{thm}{Theorem}[section]
\newtheorem*{thm*}{Theorem}

\newtheorem{prop}[thm]{Proposition}
\newtheorem{cor}[thm]{Corollary}
\newtheorem*{cor*}{Corollary}
\newtheorem*{prop*}{Proposition}

\newtheorem*{lemma*}{Lemma}
\newtheorem*{claim*}{Claim}

\theoremstyle{definition}

\newtheorem*{exmp*}{Example}
\newtheorem*{defn*}{Definition}
\newtheorem*{rem*}{Remark}
\newtheorem*{note*}{Note}

\begin{document}

\title{Complements of tori in $\#_{2k}S^2 \times S^2$ that admit a hyperbolic structure} 

\author{Hemanth Saratchandran}
\email{hemanth.saratchandran@maths.ox.ac.uk}
\date{\today}

\maketitle 
\parskip=0.2cm
\parindent=0.0cm

\begin{abstract}
We construct examples of codimension two hyperbolic link complements in closed smooth 4-manifolds with 
homeomorphism type
$\#_{2k}S^2 \times S^2$. All our examples are based on a construction of J. Ratcliffe and S. Tschantz, who
constructed 1171 non-compact finite volume hyperbolic 4-manifolds of minimal volume. We then give necessary 
conditions for a closed smooth simply connected 4-manifold to contain a codimension two link complement that
admits a hyperbolic structure. 
\end{abstract}

\tableofcontents

\section{Introduction}\label{intro}

In the early 80's Michael Freedman stunned the mathematical world by obtaining a complete classification
of all closed simply connected topological 4-manifolds (\cite{freedman}). He was able to show that the 
homeomorphism
type of a closed simply connected topological 4-manifold was characterised by two invariants. The first
one was the intersection form
\[ \langle\hspace{0.1cm} , \hspace{0.1cm}\rangle : H_2 \otimes H_2 \rightarrow \Z \]
which is a symmetric bilinear form having determinant $\pm1$. The second is the Kirby-Seibenmann invariant, which
is an element 
\[ \sigma \in H^4(M; \Z_2) \cong \Z_2 \]
that vanishes if and only if the product $M \times \R$ admits a smooth structure (see \cite{scorpan} p.221-222). 
Furthermore, he was
able to show that given an arbitrary integral symmetric bilinear form $Q$ with determinant $\pm1$,
there exists a closed simply connected topological 4-manifold that has $Q$ as its intersection form. 
The number of definite forms is enormous, coupled with the classification theorem
of Freedman this shows that the number of closed simply connected topological 4-manifolds whose intersection form
is definite is enormous. One could try and exclude some of these forms by restricting to those manifolds
that admit a smooth structure, but unfortunately Freedman's theorem does not say anything about a diffeomorphism
classification of smooth 4-manifolds, and on the existence of smooth structures. During that time there was no
way of ruling out the possibility that the enormous number of topological 4-manifolds with definite intersection
form admitted a smooth structure.  \\
Just one year after Freedman's work a second revolution occurred in the 4-manifold world. 
Simon Donaldson (\cite{donaldson}), using techniques from gauge theory, was able to show that the only definite forms
that arose as the intersection form of a closed smooth simply connected 4-manifold were 
\[ \oplus_n[+1] \text{ and} \oplus_n[-1] \]
Along with the work of Freedman, this immediately implied that the overwhelming crowd of topological 4-manifolds 
with definite intersection forms was mainly populated by ones that did not admit a smooth structure. A result
many mathematicians did not see coming. Freedman's work together with the work of Donaldson refines the 
classification theorem for closed smooth simply connected 4-manifolds, and can be expressed in the following simple
form. 
\textit{Every closed smooth simply connected 4-manifold is homeomorphic to either }
\[ S^4 \text{ or } \#_n\CP^2\#_m\overline{\CP}^2 \text{ or } \#_{\pm m}\mathcal{M}E_8\#_nS^2 \times S^2 \]
Here $\mathcal{M}E_8$ denotes the non-smoothable topological 4-manifold with the $E_8$ intersection form.
Although this theorem may seem somewhat unsatisfactory, since many of the 
$\#_{\pm m}\mathcal{M}E_8\#_nS^2 \times S^2$ are non-smoothable. It has proved to be invaluable in our understanding
of 4-manifolds.

The aim of this paper is to try and understand which homeomorphism types, in the Donaldson-Freedman classification
theorem, can contain a codimension two link complement that admits a hyperbolic structure. By hyperbolic structure
we will always mean a finite volume hyperbolic structure.
The first piece of work that was done in investigating this type of question was carried out by 
D. Ivan$\check{s}$i$\acute{c}$ in \cite{ivansic}. In that paper Ivan$\check{s}$i$\acute{c}$ constructed 
a comdimension two link complement in a closed smooth 4-manifold that was homeomorphic to $S^4$. The starting
point for his construction was to analyse the orientable double cover of a non-orientable finite volume hyperbolic 
4-manifold constructed by J. Ratcliffe and S. Tschantz. In \cite{ratcliffe} Ratcliffe and Tschantz construct
1171 non-compact finite volume hyperbolic 4-manifolds, with 22 of them being orientable and 1149 being non-
orientable, 
each of Euler characteristic $1$. Ivan$\check{s}$i$\acute{c}$ shows that the orientable double cover of the
one numbered 1011 has five cusps, with each cusp cross-section given by the 3-torus.    
By gluing in a solid 3-torus to each cusp cross-section via
the identity map, he showed that you would obtain a closed smooth simply connected 4-manifold whose Euler 
characteristic was $2$. Using the Donaldson-Freedman classification theorem he was able to conclude that the closed
smooth simply connected 4-manifold he obtained, through this gluing process, was homeomorphic to $S^4$.
In this way he had produced the first example of a codimension two hyperbolic link complement in a 4-manifold
that was homeomorphic to $S^4$. Soon after D. Ivan$\check{s}$i$\acute{c}$, J. Ratcliffe, and S. Tschantz 
constructed several more examples of codimension two hyperbolic link complements in 4-manifolds that
were homeomorphic to $S^4$ (\cite{tschantz}). The starting point of this paper is to construct many more examples
of codimension two hyperbolic link complements in closed smooth simply connected 4-manifolds, whose
homeomorphism type we can classify. We show how to construct codimension two link complements in smooth 4-manifolds 
whose homeomorphism type is given by $\#_{2k}(S^2 \times S^2)$ for $k \geq 1$.

\begin{thm}
For $k \geq 1$ there exists a collection of $8k + 5$ tori embedded in a smooth 4-manifold $X$ such that
$X$ is homeomorphic to $\#_{2k}(S^2 \times S^2)$, and $X-L$ admits a finite volume hyperbolic structure.
\end{thm}

Our approach to proving the above theorem starts by constructing certain finite covers of the orientable double
cover of the hyperbolic manifold numbered 1011 in the list constructed by Ratcliffe and Tschantz. The cusps
of each of these covers each have a cross-section given by the 3-torus. We show that if we glue in
a solid 3-torus, via the identity, we obtain a closed smooth simply connected 4-manifold that admits a spin
structure. Then using work of D. Long and A. Reid on the signature of a non-compact finite volume
hyperbolic 4-manifold (\cite{long}), we prove that the signature of these closed
4-manifolds must vanish. These two results are enough to classify the homeomorphism type of the closed
4-manifolds, obtained from gluing in solid 3-tori to each cusp cross-section, and obtain the above theorem. 

The second, and final, goal of the paper is to understand whether all the homeomorphism types of 
smooth 4-manifolds, 
occurring in the Donaldson-Freedman classification theorem, can contain a codimension two link complement
that admits a hyperbolic structure. For example, one of the homeomorphism types occurring in the
Donaldson-Freedman classification theorem is given by $\#_{-2}\mathcal{M}E_8 \#_3S^2 \times S^2$. The fact
that this admits a smooth structure follows from the observation that it has the same 
intersection form, given by $-2E_8 \oplus 3H$ (where $H$ denotes the hyperbolic form), as a $K3$ surface. One can
then ask if a $K3$ surface contains a codimension two link complement that admits a hyperbolic structure?
We are able to show that this cannot be the case, and in fact many homeomorphism types, occurring in the
Donaldson-Freedman classification, cannot contain a codimension two link complement that admits a 
hyperbolic structure. Our main theorem, in regards to this second goal, takes the following form.

\begin{thm}
Let $M$ be a closed smooth simply connected 4-manifold that has a codimension two link complement that admits a 
hyperbolic structure. Then the homeomorphism
type of $M$ falls into one of the following three categories:

\begin{itemize}
\item $S^4$

\item $\#_k(S^2 \times S^2)$, $k > 0$.

\item $\#_k\CP^2 \#_k \overline{\CP^2}$, $k > 0$.
\end{itemize}
\end{thm}

We prove this theorem by showing that any such closed smooth 4-manifold that contains a codimension two link 
complement that admits a hyperbolic structure must have vanishing signature. This immediately gives the above
result.


\section*{Acknowledgements}
I would like to thank Andras Juhasz and Marc Lackenby for finding an error in an earlier version of this work, and
for their help in repairing this error.
I would like to thank Igor Belegradek for providing a useful reference on eta invariants.
Finally, I must thank John Parker for his valuable comments on an earlier version of this work

\section{The Ratcliffe-Tschantz manifolds}\label{prelim}

In this section we will introduce the reader to the Ratcliffe-Tschantz hyperbolic 4-manifolds.
Each of their hyperbolic 4-manifolds has the hyperbolic 24-cell, a four dimensional self dual polyhedron, 
as a fundamental domain. In view of this we start by describing how the hyperbolic 24-cell is constructed, and
the notation we will be using to describe it. We then state the main theorem we will need from their paper
\cite{ratcliffe}, and explain the coding they developed to identify each of the 1171 hyperbolic 4-manifolds
they constructed.

Let $S_{(*,*,*,*)}$ denote a sphere of radius 1 centred at a point in $\R^4$ whose coordinates have two $\pm 1$'s 
and whose other two coordinates are both zero.
For example $S_{(+1,+1,0,0)}$ denotes the sphere of radius 1 centred at the point $(1,1,0,0)$ in $\R^4$. If we let 
$\mathbb{H}^4$ denote the ball model 
of hyperbolic 4-space then we find that all the spheres $S_{(*,*,*,*)}$ intersect the sphere at infinity 
orthogonally. 
This implies that each such sphere determines a hyperplane in $\mathbb{H}^4$. 
If we let $Q_{(*,*,*,*)}$ denote the corresponding half-space that contains the origin, and then take the 
intersection of all such half-spaces, we find that we obtain
a 24-sided polyhedron in $\mathbb{H}^4$. This polyhedron is known as the hyperbolic 24-cell, and we will denote 
it by $P$. It is a four dimensional self dual polyhedron.
We will denote the side of $P$ that lies on the sphere
$S_{(*,*,*,*)}$ also by $S_{(*,*,*,*)}$. All the dihedral angles of $P$ are $\pi/2$ and it has twenty-four 
ideal vertices. We can explicitly describe each ideal
vertex: We have 8 vertices of the form $v_{(\pm 1, 0,0,0)} = (\pm 1,0,0,0)$, $v_{(0, \pm 1, 0,0)} = (0,\pm 1,0,0)$, 
$v_{(0,0, \pm 1, 0)} = (0,0,\pm 1,0)$, 
$v_{(0,0,0, \pm 1)} = (0,0,0,\pm 1)$ plus 16 vertices of the form 
$v_{(\pm 1/2, \pm 1/2, \pm 1/2, \pm 1/2)} = (\pm 1/2, \pm 1/2, \pm 1/2, \pm 1/2)$. 

Let $\Gamma^n$  be the group of $(n+1) \times (n+1)$ positive Lorentzian matrices with integer entries. The group 
$\Gamma^n$ is an infinite discrete subgroup of the group $O(n, 1)$ of
Lorentzian $(n+1) \times (n+1)$ matrices. The principal congruence two subgroup of $\Gamma^n$ is the group 
$\Gamma^n_2$ of all matrices in $\Gamma^n$ that are congruent
to the identity matrix modulo two.
In their paper \cite{ratcliffe} Ratcliffe and Tschantz classify all the hyperbolic space forms $\mathbb{H}^4/\Gamma$ 
where $\Gamma$ is a torsion free subgroup
of minimal index in the group $\Gamma^4_2$ (they actually do the case $n = 2, 3$ as well but as we are only 
interested in the $n = 4$ case we will not worry about these
other two). Their main theorem takes the form:

\begin{thm}
There are, up to isometry, exactly 1171 hyperbolic space-forms $\mathbb{H}^4/\Gamma$ where $\Gamma$ is a torsion free subgroup of minimal index in the group
$\Gamma^4_2$. Only 22 of these manifolds are orientable.
\end{thm}

Each of the 1171 hyperbolic 4-manifolds coming from the above theorem have Euler characteristic $1$.
The proof of this theorem which leads to the construction of the 1171 hyperbolic space-forms $\mathbb{H}^4/\Gamma$ 
involves finding suitable side pairings of the 24-cell that will give rise to hyperbolic 4-manifolds. 
The interested reader can consult p.111 of \cite{ratcliffe} for the details of the proof.

All their side pairings are of the form $rk$, where $k$ is a diagonal matrix with diagonal entry taking the form 
$(\pm 1, \pm 1, \pm 1, \pm 1)$, which is to be interpreted as 
a composition of reflections in the coordinate planes $x_i = 0$ for $1 \leq i \leq 4$. In order to
identify a particular element $k$ it suffices to simply give the diagonal  $(\pm 1, \pm 1, \pm 1, \pm 1)$, where a 
$-1$ in some position tells us that we reflect in the coordinate 
corresponding to that position, and a $+ 1$ tell us that we do nothing. For example $(-1, +1, -1, +1)$ is the 
composition of reflection in $x_1 = 0$ followed by reflection in $x_3 = 0$.

Ratcliffe and Tschantz develop a coding system for the $k$-part of each side pairing transformation which we now 
describe. 
The polyhedron $P$ has 24 three dimensional sides, since a side pairing transformation must identify pairs of sides 
we need to give twelve transformations. We will
denote these transformations by the letters $a, b, \ldots , k, l$.
We then group the letters $a, b, \ldots ,l$ and the sides of $P$ into the following groups:
\[ \{a, b, S_{(\pm 1,\pm 1, 0,0)}\}, \{c, d, S_{(\pm 1, 0, \pm 1, 0)}\}, \{e, f, S_{(0,\pm 1, \pm 1,0)}\}, \]
\[\{g, h, S_{(\pm 1, 0, 0, \pm 1)}\}, \{i, j, S_{(0,\pm 1, 0, \pm 1)}\}, \{k, l, S_{(0, 0, \pm 1, \pm 1)}\}. \]
In each of the above sets the letters are pairings between the spheres in that set (remember the spheres represent the sides of the 24-cell $P$ by the way we constructed it). We give 
the sides the ordering $(+1, +1) < (+1, -1) < (-1, +1) < (-1, -1)$. So for example taking the first group above our order tells us that 
\[ S_{(+1, +1, 0,0)} <  S_{(+1,-1, 0,0)} < S_{(-1, +1, 0,0)} < S_{(-1, -1, 0,0)}  \]
The first letter always pairs the side with $+1, +1$ to one of the other sides (we will show how to determine this side), and the second letter pairs the next unused side
(with respect to the above ordering) with the last side.
The actual side pairing transformations are encoded by a string of six characters from the set 
\[ \{1,2,3,4,5,6,7,8,9,A,B,C,D,E,F\} \]
one character for each of the above sets. Each character represents a particular $k$-part, and the following table shows the correspondence:

\begin{center}
\begin{tabular}{|l | l | }
\hline	

\textbf{Character} &  \textbf{k-part}  \\ \hline

1 & $k_{_{(-1, +1, +1, +1)}}$  \\ \hline

2 & $k_{_{(+1, -1, +1, +1)}}$  \\ \hline

3 & $k_{_{(-1, -1, +1, +1)}}$  \\ \hline

4 & $k_{_{(+1, +1, -1, +1)}}$ \\ \hline

5 & $k_{_{(-1, +1, -1, +1)}}$ \\ \hline
 
6 & $k_{_{(+1, -1, -1, +1)}}$ \\ \hline 
 
7 & $k_{_{(-1, -1, -1, +1)}}$ \\ \hline 
 
8 & $k_{_{(+1, +1, +1, -1)}}$ \\ \hline  
 
9 & $k_{_{(-1, +1, +1, -1)}}$ \\ \hline   
 
A & $k_{_{(+1, -1, +1, -1)}}$ \\ \hline 
 
B & $k_{_{(-1, -1, +1, -1)}}$ \\ \hline  
 
C & $k_{_{(+1, +1, -1, -1)}}$ \\ \hline   
 
D & $k_{_{(-1, +1, -1, -1)}}$ \\ \hline    
 
E & $k_{_{(+1, -1, -1, -1)}}$ \\ \hline  
 
F & $k_{_{(-1, -1, -1, -1)}}$ \\ \hline

\end{tabular} 
\end{center}

The coding of a particular manifold will take the form of six characters from the above set of characters. For 
example, a code can look like
\[ \textbf{1428BD} .\]
Each character in the above code tells us what the $k$-part of each pair of transformations in the group
\[ \{(a, b), (c, d), (e, f), (g, h), (i, j), (k, l)\} \]
is. For example for the code \textbf{1428BD} the first character is \textbf{1}, this tells us that for the 
particular manifold corresponding to this code
the side pairings $a$ and $b$ have $k$-parts given by $k_{_{(-1, +1, +1, +1)}}$. The next character in the code is 
\textbf{4}, this tells us that the side pairings
$c$ and $d$ have $k$-part  $k_{_{(+1, +1, -1, +1)}}$. Continuing in this way we can determine all the $k$-parts of 
each of the side pairings $a, b, c, d, \ldots , k, l$.

We have still not explained what the $r$-part of the side pairing transformations are, they are just reflections in 
the image side, viewed as a hyperplane in $\mathbb{H}^4$.
For example for the code \textbf{1428BD}, described above, we know that the side pairing transformation $a$ pairs 
side $S_{(+1,+1,0,0)}$ to side $S_{(-1,+1,0,0)}$, and
has $k$-part given by $k_{_{(-1,+1,+1,+1)}}$. Using the fact that the $r$-part is just reflection in the image side 
we have that
\[ a := rk_{_{(-1,+1,+1,+1)}} \]
where $r$ is reflection in the side $S_{(-1,+1,0,0)}$.

Towards the end of their paper (see p. 117-124 \cite{ratcliffe}) Ratcliffe and Tschantz have tables of the 1171 
manifolds giving various information about these 
manifolds, in particular the tables tell us what the side pairing code for each manifold is, and the link structure 
of each cusp. I.e. the flat structure on each cusp cross-section. 
The notation they use to represent these link types is as follows. 
$\textbf{A}, \textbf{B}, \ldots , \textbf{J}$ represent the ten closed Euclidean 3-manifolds in the order given by 
Hantzsche and Wendt in  \cite{hantzsche}. These are also given in Wolf's text \cite{wolf} p.122, where
in his notation $\mathcal{G}_1$, \ldots , $\mathcal{G}_6$ corresponds to $\textbf{A}, \ldots , \textbf{F}$ 
respectively, and $\mathcal{B}_1$, \ldots , $\mathcal{B}_4$ corresponds
to $\textbf{G}, \ldots , \textbf{J}$. The orientable ones are given by $\textbf{A}, \ldots , \textbf{F}$ and the 
non-orientable ones are $\textbf{G}, \ldots , \textbf{J}$. $\textbf{A}$ is the 3-torus, and $\textbf{B}$ is an 
orientable 
$S^1$-fibre bundle over the Klein bottle. The non-orientable manifolds $\textbf{G}$ and $\textbf{H}$ both
have $\textbf{A}$ as their orientable double cover, $\textbf{I}$ and $\textbf{J}$ have $\textbf{B}$ as their
orientable double cover. Furthermore, $\textbf{A}$ and $\textbf{B}$ are the only orientable ones that
are $S^1$-fibre bundles over a compact surface, all the non-orientable ones are also $S^1$-fibre bundles
over a compact surface.

\section{The orientable double cover of Manifold 1011}\label{double_cover}

The Ractliffe-Tschantz manifolds are constructed by taking the hyperbolic 24-cell and applying a particular group
of side pairing transformations. Using this fact we showed in \cite{sarat} how to construct a Kirby diagram for any 
one of the Ratcliffe-Tschantz manifolds. The basic
idea of the construction was that each equivalence class of a codimension $k$ side, under the side pairing
tranformation group, corresponds to a $k$-handle. 
In this section we are going to briefly explain how to construct a Kirby diagram for the orientable double cover of 
the non-orientable manifold numbered 1011 in the Ratcliffe-Tschantz census. In the section that follows we will
use the Kirby diagram to show that the orientable double cover of manifold 1011 must be spin.

The starting point to obtaining a Kirby diagram is to observe that
each side of the polyhedron $P$ corresponds to a sphere of the form $S_{(*,*,*,*)}$, which itself is identified by 
the co-ordinate $(*,*,*,*)$. Radially projecting this point to the boundary $\partial{B^4} = S^3$, we obtain
a collection of points on $S^3$ corresponding to each side of the hyperbolic 24-cell $P$. 
We define the map:
\[ \phi : S^3 \rightarrow \R^4 \cup \{\infty\} \]
by 
\[ \phi(x_1, x_2, x_3, x_4) = (0, 0, 0, 1) + \frac{2}{x_1^2 + x_2^2 + x_3^2 + (x_4 - 1)^2}(x_1, x_2, x_3, x_4 - 1) 
\]
Using the above map we can map each point in $S^3$ that corresponds to a sphere onto a point in $\R^3$. For example 
if we take the sphere $S_{(+1,0,+1,0)}$, then we know that
the associated point on $S^3$ has coordinates $(1/\sqrt{2}, 0, 1/\sqrt{2}, 0)$. Applying $\phi$ to this point we get 
the point  $(1/\sqrt{2}, 0, 1/\sqrt{2})$. Doing this
for all the points corresponding to the sides of $P$ we get a collection of points in $\R^3$, each one corresponding
to a side of the hyperbolic 24-cell $P$. As the polyhedron $P$ is self-dual the map $\phi$ gives us a way to 
visualise the boundary of $P$ in $\R^3$. The collection of points in $\R^3$, obtained from the map $\phi$, 
correspond to attaching spheres of a Kirby diagram. The attaching regions of the 2-handles can be found
by seeing how each codimension 2 side maps to $\R^3$ under $\phi$ (see \cite{sarat} for details).

The hyperbolic 4-manifold we will be focusing on is numbered 1011 in the Ratcliffe-Tschantz census. 
From here on in we will denote manifold 1011 by $M$, and its orientable double cover by $\widetilde{M}$.

The side pairing code for $M$ is \textbf{14FF28}, and the explicit side pairings are given as:

\[\xymatrixcolsep{5pc} \xymatrix{ S_{(+1,+1,0,0)}  \ar[r]^a_{k_{(-1,+1,+1,+1)}} & S_{(-1,+1,0,0)} } \hspace{2cm} \xymatrix{S_{(+1,-1,0,0)}  \ar[r]^b_{k_{(-1,+1,+1,+1)}} & S_{(-1,-1,0,0)} } \]

\[\xymatrixcolsep{5pc} \xymatrix{ S_{(+1,0,+1,0)}  \ar[r]^c_{k_{(+1,+1,-1,+1)}} & S_{(+1,0,-1,0)} } \hspace{2cm} \xymatrix{S_{(-1,0,+1,0)}  \ar[r]^d_{k_{(+1,+1,-1,+1)}} & S_{(-1,0,-1,0)} } \]

\[\xymatrixcolsep{5pc} \xymatrix{ S_{(0,+1,+1,0)}  \ar[r]^e_{k_{(-1,-1,-1,-1)}} & S_{(0,-1,-1,0)} } \hspace{2cm} \xymatrix{S_{(0,+1,-1,0)}  \ar[r]^f_{k_{(-1,-1,-1,-1)}} & S_{(0,-1,+1,0)} } \]

\[\xymatrixcolsep{5pc} \xymatrix{ S_{(+1,0,0,+1)}  \ar[r]^g_{k_{(-1,-1,-1,-1)}} & S_{(-1,0,0,-1)} } \hspace{2cm} \xymatrix{S_{(+1,0,0,-1)}  \ar[r]^h_{k_{(-1,-1,-1,-1)}} & S_{(-1,0,0,+1)} } \]

\[\xymatrixcolsep{5pc} \xymatrix{ S_{(0,+1,0,+1)}  \ar[r]^i_{k_{(+1,-1,+1,+1)}} & S_{(0,-1,0,+1)} } \hspace{2cm} \xymatrix{S_{(0,+1,0,-1)}  \ar[r]^j_{k_{(+1,-1,+1,+1)}} & S_{(0,-1,0,-1)} } \]

\[\xymatrixcolsep{5pc} \xymatrix{ S_{(0,0,+1,+1)}  \ar[r]^k_{k_{(+1,+1,+1,-1)}} & S_{(0,0,+1,-1)} } \hspace{2cm} \xymatrix{S_{(0,0,-1,+1)}  \ar[r]^l_{k_{(+1,+1,+1,-1)}} & S_{(0,0,-1,-1)} }. \]

Using the map $\phi$ we can map each of the
sides of $M$ to $\R^3$. These points in $\R^3$ correspond to the 1-handles of $M$.
So as to make our lives easier when dealing with $\widetilde{M}$ we are going to label
the above sides of the 24-cell by capital letters. The following table shows the labellings
we choose, and the associated points in $\R^3$ corresponding to each side under the map $\phi$.

\begin{tabular}{|l|l|l||l|l|l|}
$A$ &  $S_{++00}$ & $(\frac{1}{\sqrt{2}}, \frac{1}{\sqrt{2}}, 0)$  &  $A'$ & $S_{-+00}$ & $(\frac{-1}{\sqrt{2}}, 
\frac{1}{\sqrt{2}}, 0)$  \\
$B$ &  $S_{+-00}$ & $(\frac{1}{\sqrt{2}}, \frac{-1}{\sqrt{2}}, 0)$ &  $B'$ & $S_{--00}$  & $(\frac{-1}{\sqrt{2}}, 
\frac{-1}{\sqrt{2}}, 0)$ \\
$C$ &  $S_{+0+0}$ &  $(\frac{1}{\sqrt{2}}, 0, \frac{1}{\sqrt{2}})$ &  $C'$ & $S_{+0-0}$ & $(\frac{1}{\sqrt{2}}, 0, 
\frac{-1}{\sqrt{2}})$ \\
$D$ &  $S_{-0+0}$ &  $(\frac{-1}{\sqrt{2}}, 0, \frac{1}{\sqrt{2}})$ &  $D'$ & $S_{-0-0}$ & $(\frac{-1}{\sqrt{2}}, 
0, \frac{-1}{\sqrt{2}})$ \\ 
$E$ &  $S_{0++0}$ &  $(0, \frac{1}{\sqrt{2}} ,\frac{1}{\sqrt{2}})$ &  $E'$ & $S_{0--0}$ &  $(0, \frac{-1}
{\sqrt{2}} 
,\frac{-1}{\sqrt{2}})$ \\ 
$F$ &  $S_{0+-0}$ &   $(0, \frac{1}{\sqrt{2}} ,\frac{-1}{\sqrt{2}})$ &  $F'$ & $S_{0-+0}$ &  $(0, \frac{-1}
{\sqrt{2}} ,\frac{1}{\sqrt{2}})$ \\ 
$G$ &  $S_{+00+}$ &  $(1 + \sqrt{2}, 0, 0)$ &  $G'$ & $S_{-00-}$ & $(1 - \sqrt{2}, 0, 0)$ \\ 
$H$ &  $S_{+00-}$ &  $(-1 + \sqrt{2}, 0, 0)$ &  $H'$ & $S_{-00+}$ & $(-1 - \sqrt{2}, 0, 0)$ \\
$I$ &  $S_{0+0+}$ &  $(0, 1 + \sqrt{2}, 0)$ &  $I'$ & $S_{0-0+}$ &  $(0, -1 - \sqrt{2}, 0)$ \\
$J$ &  $S_{0+0-}$ &   $(0, -1 + \sqrt{2}, 0)$ &  $J'$ & $S_{0-0-}$ &  $(0, 1 - \sqrt{2}, 0)$  \\
$K$ &  $S_{00++}$ &   $(0, 0, 1 + \sqrt{2})$ &  $K'$ & $S_{00+-}$ &  $(0, 0, -1 + \sqrt{2})$   \\
$L$ &  $S_{00-+}$ &   $(0, 0, -1 - \sqrt{2})$  &  $L'$ & $S_{00--}$ &  $(0, 0, 1 - \sqrt{2})$  
\end{tabular}

We are not going to show how the Kirby diagram of $M$ looks like as we do not need it. The interested reader
can consult \cite{sarat} p.20, where we give details on how to construct the Kirby diagram for $M$.

The manifold has five cusps given by the code \textbf{GGGGG}, see end of section \ref{prelim} for the definition
of this notation (or in \textit{Wolf's} notation 
$\mathcal{B}_1 \mathcal{B}_1 \mathcal{B}_1 \mathcal{B}_1 \mathcal{B}_1$). 
In \cite{sarat_2} p.16 we showed
how to work out the parabolic transformations corresponding to each such cusp. The following table
shows the generators of each of the parabolic subgroups associated to each of the five cusps. 

\begin{tabular}{|l | l | }
\hline	

Ideal vertex & Generators of parabolic subgroup \\ \hline

$\{(1,0,0,0), (-1,0,0,0)\}$ & $c$, $a^{-1}b$, $a^{-1}g$ \\ \hline

$\{(0,1,0,0), (0,-1,0,0)\}$ & $a$, $e^{-1}f$, $e^{-1}i$ \\ \hline

$\{(0,0,1,0), (0,0,-1,0)\}$ & $k$, $c^{-1}d$, $c^{-1}e$ \\ \hline

$\{(0,0,0,1), (0,0,0,-1)\}$ & $j$, $g^{-1}h^{-1}$, $g^{-1}k$ \\ \hline

 $\{(\pm 1/2, \pm 1/2, \pm 1/2, \pm 1/2)\}$ & $e^{-1}g$, $a^{-1}k^{-1}ak$, $a^{-1}k^{-1}j^{-1}f$ \\ \hline

\end{tabular} \\ \\ \\

The closed flat 3-manifold
given by \textbf{G} is a non-orientable $S^1$-fibre bundle. The parabolic translations that correspond
to the $S^1$-fibre are given by $c$, $a$, $k$, $j$, and $e^{-1}g$.

The starting point
to obtaining a Kirby diagram for the orientable double cover is to identify the orientation preserving isometries and those 
that are orientation
reversing. 
Recall, any side pairing transformation is written as the composition $rk$, where $r$ is reflection in the image 
side and $k$ is a diagonal matrix.
It is clear that $r$ is orientation reversing (since it is a reflection in a hyperplane), therefore we can conclude 
that the orientation
preserving isometries are those whose $k$-part has an odd number of $-1$'s, and the orientation reversing isometries 
are those that
have an even number of $-1$'s. We then find that the
orientation preserving isometries are given by $a,b,c,d,i,k$, and the orientation reversing isometries are given 
by $e,f,g,h$.

A fundamental domain for the orientable double cover consists of two copies of the 24-cell $P$ 
attached along a codimension one side corresponding to one of the orientation reversing isometries. If we take the 
orientation reversing isometry
$g^{-1}$, then we can think of the fundamental domain of the orientable double cover to consist of the union $P \cup 
g^{-1}\cdot P$. 
In our paper \cite{sarat_3} we showed how to compute the side pairing transformations
associated to the orientable double cover of any one of the non-orientable Ratcliffe-Tschantz manifolds using the
Reidemeister-Schreier rewriting process. The reader who is interested in the details of how one does this
should take a look at that paper. 
 
The side pairing transformations for $\widetilde{M}$ are:

\[\xymatrixcolsep{4pc} \xymatrix{ A  \ar[r]^a & A' } \hspace{0.5cm} \xymatrix{g^{-1}A  \ar[r]^{g^{-1}ag} & g^{-1}A' } \hspace{2cm} \xymatrix{ B  \ar[r]^b & B' }  \hspace{0.5cm} 
 \xymatrix{g^{-1}B  \ar[r]^{g^{-1}bg} & g^{-1}B' } \]

\[\xymatrixcolsep{4pc} \xymatrix{ C  \ar[r]^c & C' } \hspace{0.5cm} \xymatrix{g^{-1}C  \ar[r]^{g^{-1}cg} & g^{-1}C' } \hspace{2cm} \xymatrix{ D  \ar[r]^d & D' }  \hspace{0.5cm} 
 \xymatrix{g^{-1}D  \ar[r]^{g^{-1}dg} & g^{-1}D' } \]

\[\xymatrixcolsep{4pc} \xymatrix{ E  \ar[r]^{g^{-1}e} & g^{-1}E' } \hspace{0.5cm} \xymatrix{g^{-1}E  \ar[r]^{eg} & E' } \hspace{2cm} \xymatrix{ F  \ar[r]^{g^{-1}f} & g^{-1}F' }  \hspace{0.5cm} 
 \xymatrix{g^{-1}F  \ar[r]^{fg} & F' } \]

\[\xymatrixcolsep{4pc} \xymatrix{ G  \ar[r]^{g^{-1}g} & g^{-1}G' } \hspace{0.5cm} \xymatrix{g^{-1}G  \ar[r]^{gg} & G' } \hspace{2cm} \xymatrix{ H  \ar[r]^{g^{-1}h} & g^{-1}H' }  \hspace{0.5cm} 
 \xymatrix{g^{-1}H  \ar[r]^{hg} & H' } \]

\[\xymatrixcolsep{4pc} \xymatrix{ I  \ar[r]^i & I' } \hspace{0.8cm} \xymatrix{g^{-1}I  \ar[r]^{g^{-1}ig} & g^{-1}I' } \hspace{2.3cm} \xymatrix{J  \ar[r]^j & J'}  \hspace{0.5cm} 
 \xymatrix{g^{-1}J  \ar[r]^{g^{-1}jg} & g^{-1}J' } \]

\[\xymatrixcolsep{4pc} \xymatrix{ K  \ar[r]^k & K' } \hspace{0.5cm} \xymatrix{g^{-1}K  \ar[r]^{g^{-1}kg} & g^{-1}K' } \hspace{2cm} \xymatrix{ L  \ar[r]^l & L' }  \hspace{0.5cm} 
 \xymatrix{g^{-1}L  \ar[r]^{g^{-1}lg} & g^{-1}L' } \]

The general construction of the Kirby diagram of a Ratcliffe-Tschantz manifold from the side pairing 
tranformations was explained in detail in \cite{sarat}. The same construction outlined in that paper
also works in this case. As we will not need the details of the construction we will simply show the
reader what the Kirby diagram looks like. 

The totality of the Kirby diagram will be split into
four separate diagrams. Each diagram contains two sub diagrams corresponding to the fact that
a fundamental domain for $\widetilde{M}$ consists of two components, $P$ and $g^{-1}P$. The first diagram
shows the $x-y$ plane of the Kirby diagram coming from $P$ and $g^{-1}P$. The second and third diagrams
show the $x-z$ and $y-z$ planes respectively of the Kirby diagram coming from $P$ and $g^{-1}P$. Finally, there
are always twelve 2-handles that do not all lie in any one of these three planes, the fourth diagram shows
how these 2-handles look like. In total there are twenty four 1-handles and forty eight 2-handles, with
twelve 2-handles in each of the $x-y$, $x-z$, $y-z$ planes, and another twelve that do not all lie in any
one of these planes. 
In the following diagrams for each side $S$ we denote the component of the 1-handle corresponding to $g^{-1}S$ by 
$S-$. The attaching maps for the attaching regions of the 1-handles are given by the above side pairing
transformations. As the Kirby diagram is split into four separate diagrams we colour code the 
2-handles in each of the four diagrams by twelve distinct colours.

The following diagram shows that part of the Kirby diagram contained in the $x-y$ plane.

\centerline{\graphicspath{ {double_cover/Manifold_1011_double_cover/2-handle_cycles/} } \includegraphics[width=14cm, height=9cm]{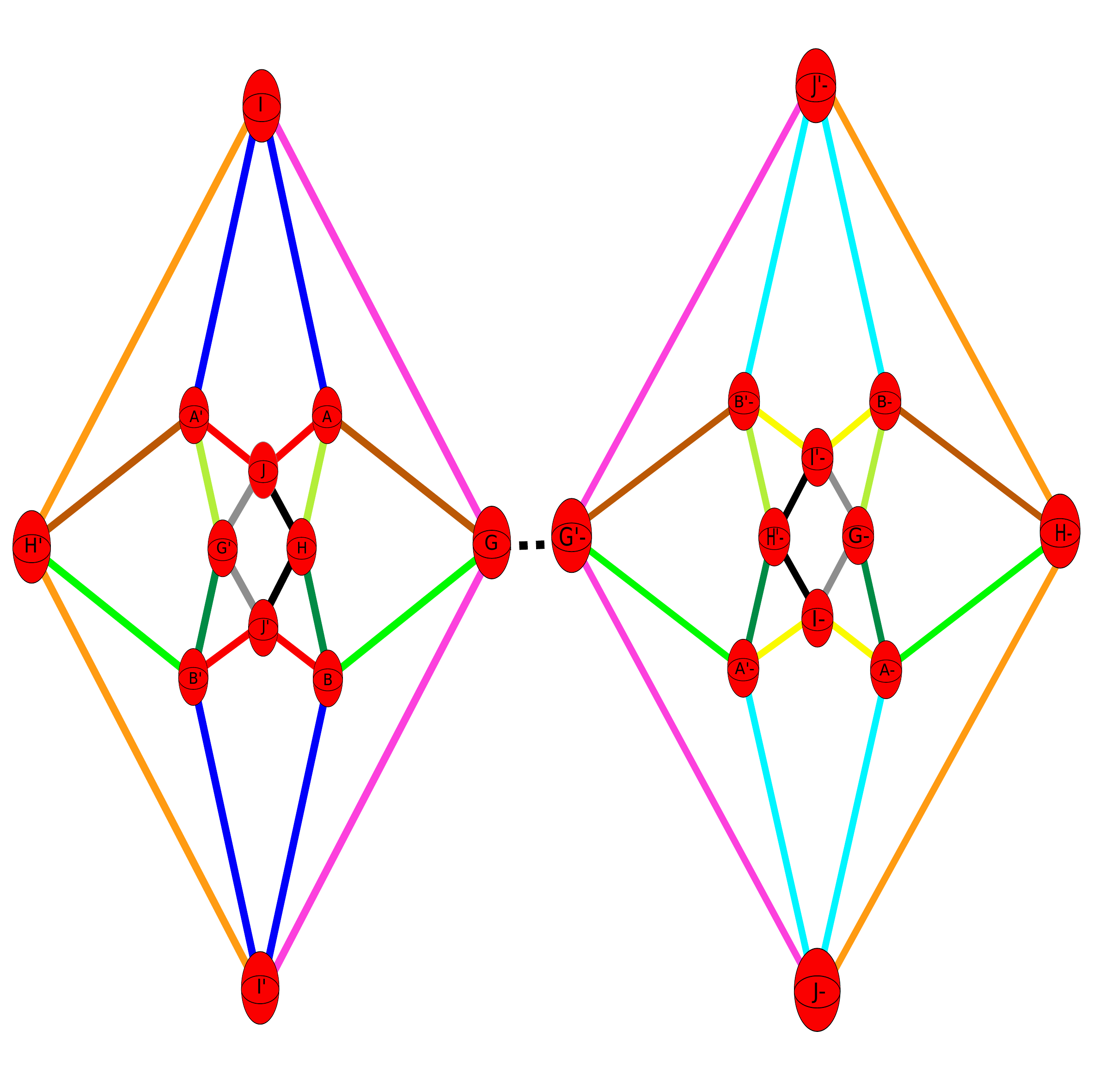}}

That part of the Kirby diagram contained in the $x-z$ plane can be seen in the following diagram.

\centerline{\graphicspath{ {double_cover/Manifold_1011_double_cover/2-handle_cycles/} } \includegraphics[width=14cm, height=9cm]{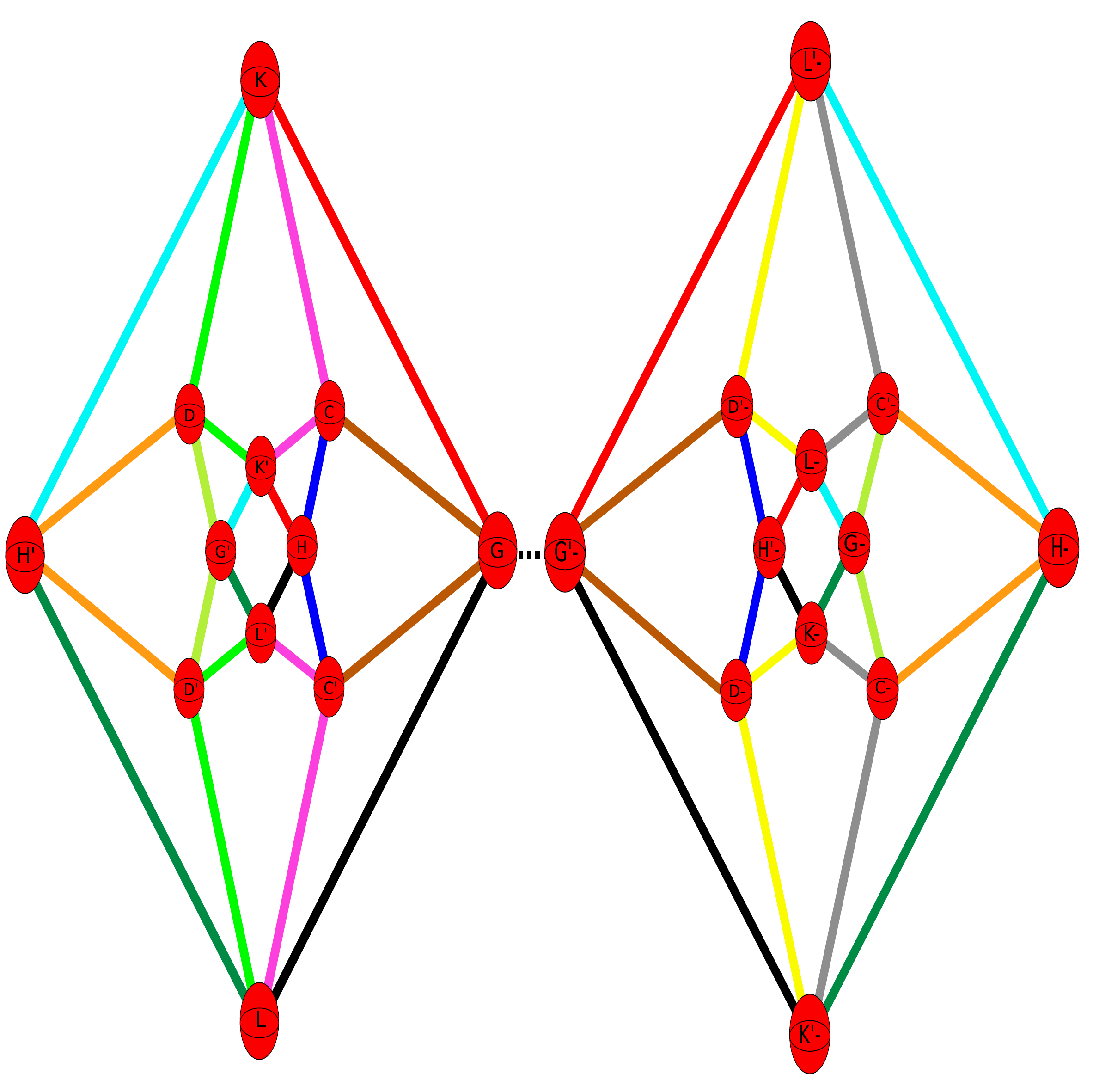}}

That part of the Kirby diagram contained in the $y-z$ plane can be seen in the following diagram.

\centerline{\graphicspath{ {double_cover/Manifold_1011_double_cover/2-handle_cycles/} } \includegraphics[width=14cm, height=9cm]{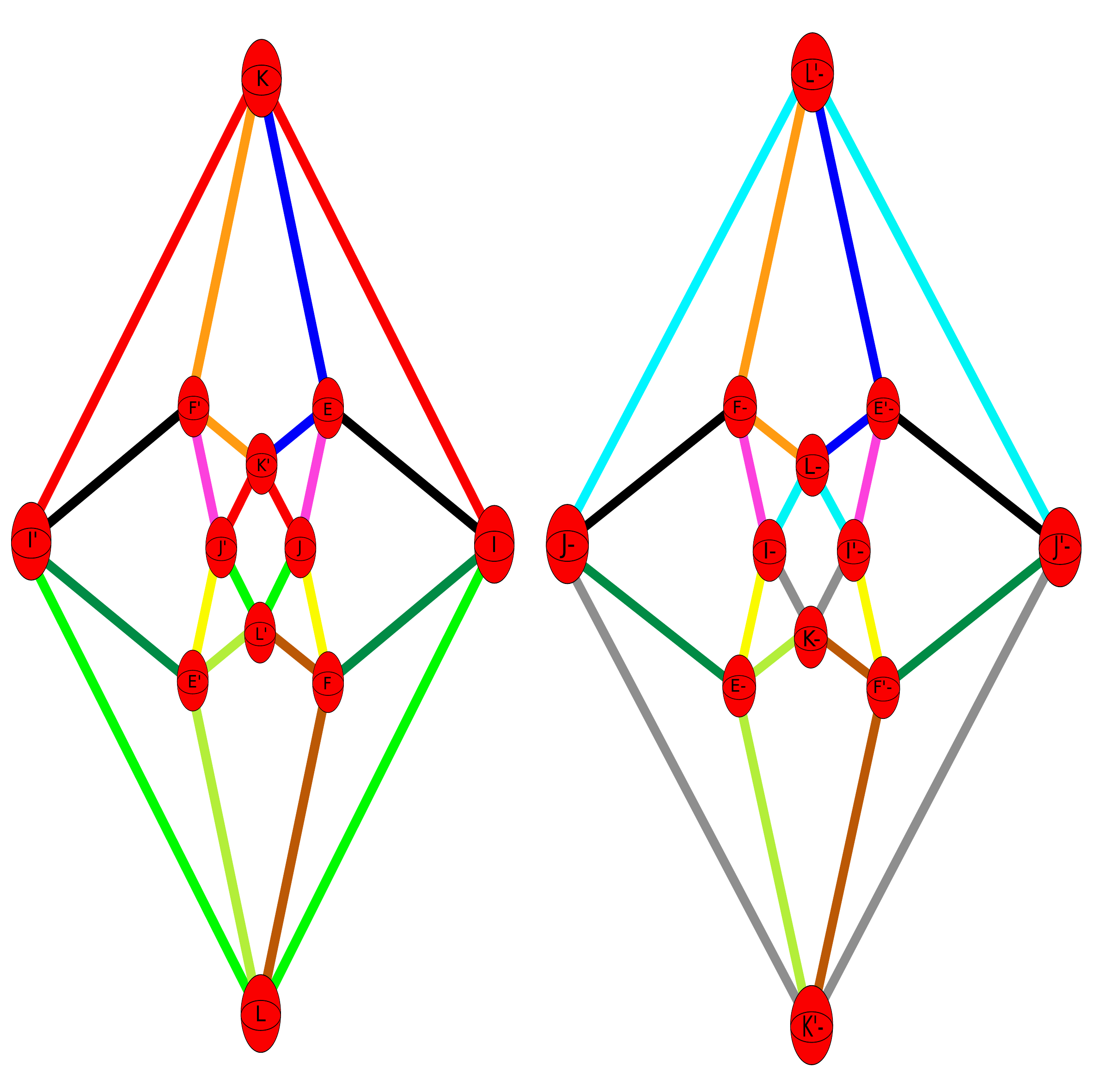}}

We said that the above shows that part of the Kirby diagram
contained in the $y-z$ plane. This is not exactly true, what we mean by this is that
in the above diagram the component on the left shows the $y-z$ plane
of the Kirby diagram coming from the $P$ part of the fundamental domain $P \cup g^{-1}P$. The component
on the right shows the $y-z$ plane of the Kirby diagram coming from the $g^{-1}P$ part.

Finally, we have the 2-handles that do not all lie in any one of the above planes, there are twelve in total.

\centerline{\graphicspath{ {double_cover/Manifold_1011_double_cover/2-handle_cycles/} } \includegraphics[width=14cm, 
height=9cm]{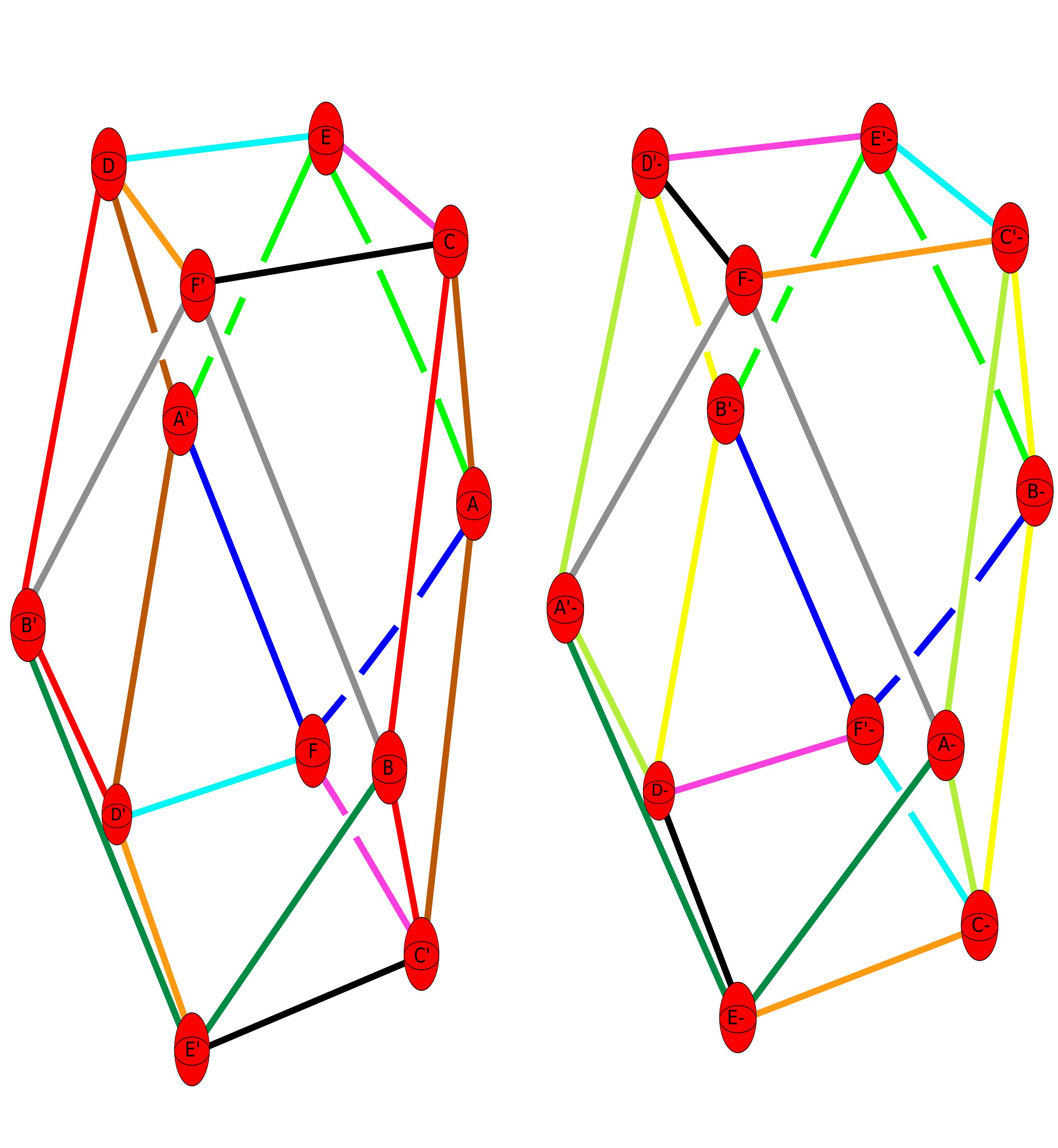}} 

Observe that in the diagrams showing the $x-y$ and $y-z$ planes, there is a 2-handle running between
$G$ and $G'-$. This 2-handle comes about due to the fact that when we form the fundamental domain
for $\widetilde{M}$ we take $P$ and transform it across the side $G$ by applying $g^{-1}$. Hence 
the 2-handle is needed to cancel the 1-handle $G, G'-$.

The double cover $\widetilde{M}$ is also a five cusped hyperbolic 4-manifold with each cusp having
type $\textbf{A}$, the 3-torus (this is because the orientable double cover of $\textbf{G}$ is $\textbf{A}$).
In \cite{ivansic} D. Ivan$\check{s}$i$\acute{c}$ proves that $\widetilde{M}$ is a codimension two hyperbolic link 
complement
in a closed smooth 4-manifold that is homeomorphic to $S^4$. 
The basic idea behind his proof is as follows.
The hyperbolic 4-manifold $M$ has five cusps each having type $\textbf{G}$, which is
an $S^1$-fibre bundle. Being an $S^1$-fibre bundle means $\textbf{G}$ bounds a natural
4-manifold, the associated disk bundle. By gluing in these disk bundles to each cusp cross-section of $M$, via the
identity, 
one obtains a closed 4-manifold for which $M$ is a hyperbolic link complement in. 
Algebraically the gluing in of disk bundles corresponds to adding certain parabolic transformations to a
presentation of the fundamental group of $M$. 
Therefore in each parabolic subgroup corresponding to each cusp there
exists a parabolic translation corresponding to the $S^1$-fibre. The parabolic translations are given by
$c$, $a$, $k$, $j$, and $e^{-1}g$. He then finds a finite presentation of the fundamental group of $M$ and
adds the relations $c = a = k = j = e^{-1}g$. He then shows that by adding these extra relations to
the fundamental group one obtains the group $\Z_2$. In other words gluing in the associated disk bundles 
to each cusp cross-section produces a closed 4-manifold with fundamental group $\Z_2$, which we denote by $N$.
It is easy to see that $N$ must have Euler characteristic equal to that of $M$, which is $1$ (remember all the 
Ratcliffe-Tschantz manifolds have Euler characteristic $1$). This is because
a flat manifold has Euler characteristic zero, so gluing in a disk bundle associated to a flat
$S^1$-fibre bundle cannot change the Euler characteristic. The orientable double cover of $N$, which we 
denote by $\widetilde{N}$, is a closed simply connected 4-manifold of Euler characteristic $2$. It is easy to see 
that $\widetilde{M}$ is a hyperbolic link complement in $\widetilde{N}$. Applying the Donaldson-Freedman 
classification theorem (see theorem \ref{simple_connec}) Ivan$\check{s}$i$\acute{c}$ is able to conclude that 
$\widetilde{N}$ is homeomorphic to $S^4$. Ivan$\check{s}$i$\acute{c}$'s theorem shows us that if we glue
in solid 3-tori to each cusp cross-section of the orientable double cover $\widetilde{M}$ we produce
a closed 4-manifold that is homeomorphic to $S^4$. We should point out that because
the 3-torus fibres over the 2-torus in more than one way the gluing in of a solid 3-torus to a cusp cross-section is 
not 
unique. It depends on a choice of what we identify as the $S^1$-fibre, or put another way the gluing in of a solid 
3-torus depends on a choice of a meridian. 

Our construction of codimension two hyperbolic link complements in 4-manifolds with homeomorphism type
$\#_{2k}(S^2 \times S^2)$ proceeds in a similar way. We construct certain finite covers of 
$\widetilde{M}$, which must have cusp cross-sections given by the 3-torus $\textbf{A}$. We then 
``fill in'' these manifolds by gluing in solid 3-tori via the identity map, with respect to a choice of meridian, to 
produce a closed 4-manifold. From here on in
we will refer to the gluing in of disk bundles to the cusp cross-sections of a hyperbolic 4-manifold via the 
identity map as a ``filling'' of that manifold.
We then prove that these fillings must be spin and have vanishing signature. An appeal to the Donaldson-Freedman
classification theorem allows us to determine the homeomorphism type of these fillings.

\section{Spin structures via the Kirby diagram}\label{chap4}\label{spin_str}

In this section we are going to prove that $\widetilde{M}$ is a spin 4-manifold. We do this
by using the construction of a Kirby diagram of $\widetilde{M}$ outlined in the previous section, and
the fact that the second Stiefel-Whitney class can be interpreted as a cocycle whose value on each 2-handle 
is the framing coefficient of that 2-handle. 

Let $E$ be a real vector bundle over a 4-manifold $X$, which we assume from here on in is orientable. We also assume 
that the fibres of $E$ have dimension $m \geq 3$, 
this will always be the case of interest for us. In the general case one can sum with a trivial bundle to obtain 
this condition. To construct a \textit{spin structure} we 
begin with a trivialisation of $E\vert_{X_1}$, where $X_1$ denotes the 1-skeleton of $X$. We wish to extend this 
trivialisation $\tau$ over each 2-handle 
$h$ of $X$. The 2-handle $h$ is a copy of $D^2 \times D^2$, which being contractible implies $E\vert_{h}$ is 
trivial, so $\tau$ over the attaching circle
determines an element of $\pi_1(SO(m)) \cong \Z_2$ (as we are assuming $m \geq 3$). It is then clear that $\tau$ 
extends over $h$ if and only if this element
in $\Z_2$ vanishes. Applying this to each 2-handle $h$ of $X$ we obtain an element of $\Z_2$ assigned to each 2-
handle. In other words, using the handle complex we
obtain a cochain $c(\tau) \in C^2(X;\Z_2)$. One can then prove that if we were to take another trivialisation 
$\tau_1$ to start with, the cochains
$c(\tau)$ and $c(\tau_1)$ will differ by a coboundary. Furthermore one can prove that for any trivialisation $\tau$, 
$c(\tau)$ is closed. This implies
we get a well-defined cohomology class denoted $w_2(E) = [c(\tau)] \in H^2(X; \Z_2)$. We can then conclude

\begin{prop}
An oriented vector bundle $E$ over $X$ admits a spin structure if and only if $w_2(E) = 0 \in H^2(X; \Z_2)$.
\end{prop}

The cohomology class $w_2(E)$ is known as the \textit{second Stiefel-Whitney class} of $E$,  
the details of why $c(\tau)$ defines a cohomology class can be found in \cite{gompf} p.180. Given an oriented 4-
manifold
$X$ its tangent bundle is an oriented vector bundle over $X$ of rank four. The second Stiefel-Whitney class of the 
tangent bundle associated to $X$ will be denoted
by $w_2(X)$.

Given a Kirby diagram of the oriented 4-manifold $X$ the following proposition shows how we can compute the second Stiefel-Whitney class $w_2(X)$.

\begin{prop}\label{spin}
For an oriented manifold $X$ given by a Kirby diagram in dotted circle notation, $w_2(X) \in H^2(X;\Z_2)$ is 
represented by the cocycle $c \in C^2(X; \Z_2)$ whose value on each 2-handle
$h$ is the framing coefficient of $h$ modulo 2.
\end{prop}

The proof of this proposition can be found in \cite{gompf} Proposition 5.7.1 and Corollary 5.7.3, p.185-186.  
The reader who is unfamiliar with dotted circle notation can
consult \cite{gompf} chap.5.6, p.167. We will just mention that one way to think of how the dotted circle method 
comes in to play is to recall that
we can cancel any 1-handle in a Kirby diagram by adding an appropriate 2-handle. Thus, adding a 1-handle to a Kirby 
diagram $X$ is the same as removing
the cancelling 2-handle. The point is that the cocore of the 2-handle corresponds to an unknotted 2-disk in $X$, 
obtained from a 2-disk in the boundary
$\partial{X}$ by pushing the interior of this 2-disk into the interior of $X$. In this way we can think of the 
addition of a 1-handle as being done by
pushing the interior of a 2-disk, originally in $\partial{X}$, into the interior of $X$ and then removing a tubular 
neighbourhood of the disk. This
tubular neighbourhood is drawn as a circle decorated with a dot, hence the name the dotted circle notation. Using 
this viewpoint of a 1-handle
one can transfer a usual Kirby diagram into dotted circle notation. On the level of the actual diagram one can 
think of the passage to dotted circle
notation as being carried out by choosing a collection of reference arcs, one for each pair of 1-handles running 
between the attaching spheres of
the 1-handle, then taking a dotted meridian to each such reference arc and collapsing the associated 1-handle via 
the reference arc. This will cause
the 2-handles running between the 1-handles to then run through the associated dotted meridian parallel to the 
choice of reference arc
(see \cite{gompf} p.169). The problem here is that there is no canonical way to do this transformation, one has to 
make a choice of how the 1-handles are going to come together to form dotted circles i.e. one has to make a choice 
of reference arc.
(see \cite{gompf} p.169 for an example of this choice). Once one fixes a choice of such reference arcs then in 
dotted circle notation all the 2-handles are knots in 3-space. They then
have a notion of a framing integer (or as it is sometimes called framing coefficient). 

The goal of the rest of this section is to use proposition \ref{spin} and show that the orientable double
cover $\widetilde{M}$ has vanishing second Stiefel-Whitney class. In order to use proposition \ref{spin}
we need to show how to transform the Kirby diagram of $\widetilde{M}$ into
dotted circle notation.

We start by showing how we can pick a system of reference arcs for each 1-handle such that when we change to dotted 
circle notation the framing of each
2-handle will be congruent to zero mod two. In order to do this we need a way to measure the framing from the 
diagram. A convenient choice to use is the \textbf{blackboard framing}, this involves projecting the knot onto a 2-
plane giving a knot diagram. The framing coefficient
is then computed as the signed number of self-crossings of the knot.
Start with the $x-y$ plane, the following diagram shows a system of reference arcs for some of the 1-handles, they 
all lie in the x-y plane and have framing zero.

\centerline {\graphicspath{ {double_cover/Manifold_1011_double_cover/framing/x-y_plane/} } \includegraphics[width=12cm, height=8cm]{reference_arcs_1}}

This leaves us with specifying reference arcs for the 1-handles $G,G'-$, $G',G-$, $H,H'-$ and $H',H-$. The reference 
arcs for some of these will
have to leave the $x-y$ plane slightly as some other handles are in the way. The following diagram shows these 
reference arcs together with
the ones shown above.

\centerline {\graphicspath{ {double_cover/Manifold_1011_double_cover/framing/x-y_plane/} } \includegraphics[width=12cm, height=9cm]{reference_arcs_2}}

Once we have chosen a particular reference arc for a 1-handle in one diagram we have to use that same reference arc for any other diagram that the
1-handle resides in. For example since the 1-handles $G,G'-$, $G',G-$, $H,H'-$ and $H',H-$ also lie in the $x-z$ plane when we choose reference arcs
for 1-handles in the $x-z$ plane we have to keep in mind that the reference arcs for $G,G'-$, $G',G-$, $H,H'-$ and $H',H-$ have already been chosen.

The following diagram shows the reference arcs for the 1-handles in the $x-z$ plane.

\centerline {\graphicspath{ {double_cover/Manifold_1011_double_cover/framing/x-z_plane/} } \includegraphics[width=12cm, height=9cm]{reference_arcs}}

The arcs you see, in the above diagram, running between the 1-handles $G,G'-$, $G',G-$, $H,H'-$ and $H',H-$ represent projections, of the
reference arcs chosen previously for these 1-handles, onto the $x-z$ plane.
The next reference arcs we want to show are for those 1-handles in the diagram corresponding to the twelve 2-handles not all lying in a single 2-plane.
To start with we show those reference arcs that have already been chosen when we considered the $x-y$ and $x-z$ planes.

\centerline {\graphicspath{ {double_cover/Manifold_1011_double_cover/framing/x-y-z_plane/} } \includegraphics[width=12cm, height=9cm]{reference_arcs_1}} 

We are left with choosing arcs for $E,E'-$, $E',E-$, $F,F'-$ and $F',F-$. For the 1-handles $E,E'-$ and $F,F'-$ we choose horizontal arcs that
run behind the diagram as the following shows.

\centerline {\graphicspath{ {double_cover/Manifold_1011_double_cover/framing/x-y-z_plane/} } \includegraphics[width=12cm, height=10cm]{reference_arcs_2}} 

For the 1-handles $E',E-$ and $F',F-$ we choose horizontal arcs running in the front as the following shows.

\centerline {\graphicspath{ {double_cover/Manifold_1011_double_cover/framing/x-y-z_plane/} } \includegraphics[width=12cm, height=10cm]{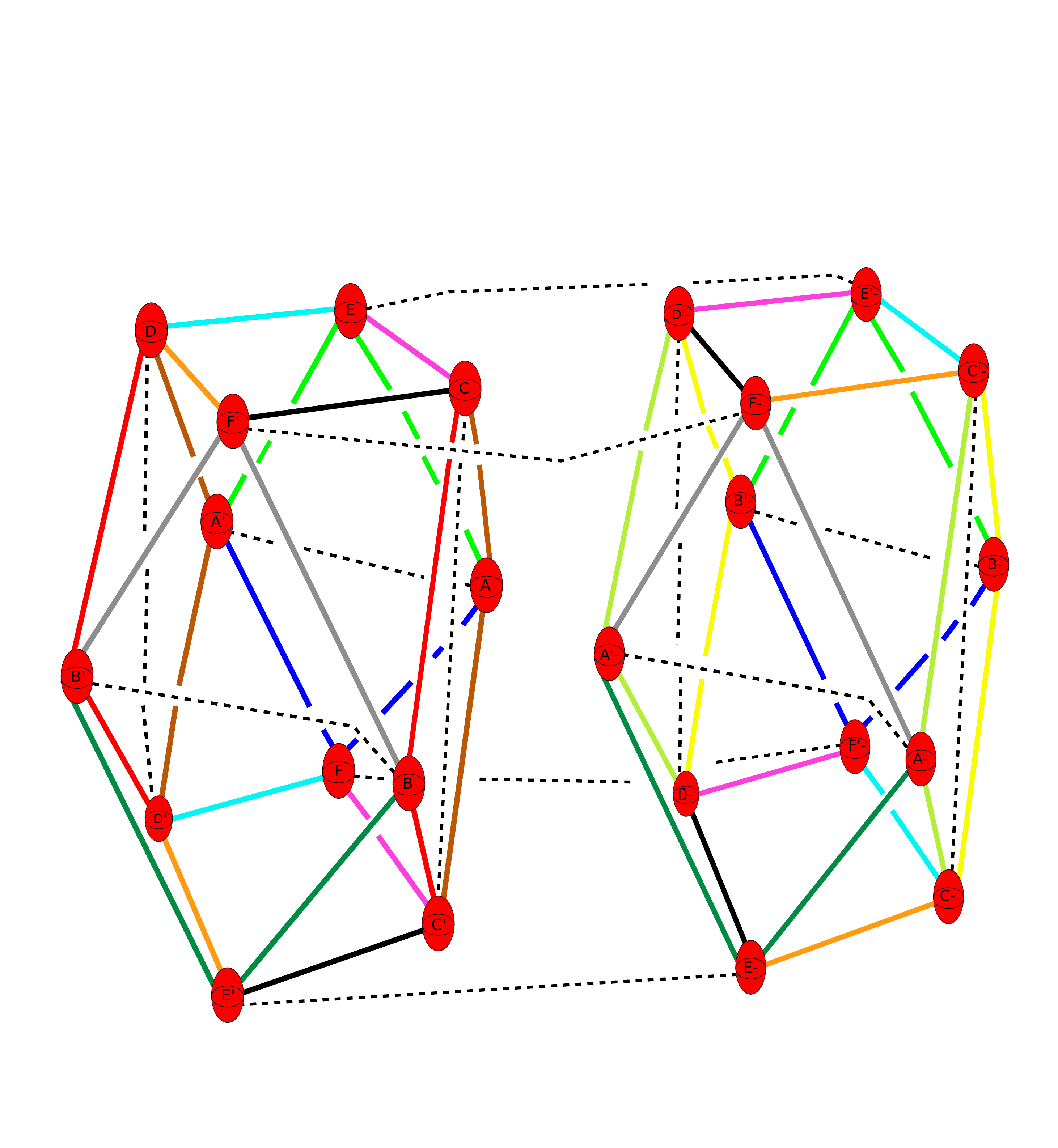}} 

The final case to consider is the diagram that shows the $y-z$ plane and the plane parallel to the $y-z$ plane. For this diagram we choose the following
arcs.

\centerline {\graphicspath{ {double_cover/Manifold_1011_double_cover/framing/y-z_plane/} } \includegraphics[width=11cm, height=8cm]{reference_arcs}} 

We have not shown how the arcs running between the 1-handles $E,E'-$, $E',E-$, $F,F'-$ and $F',F-$ look like from 
this plane. From our choice 
above it should be clear to the reader that they run horizontally outside the diagram.

We claim that with this system of reference arcs when we pass to dotted circle notation each 2-handle 
will have even framing.
This is not too hard to see, one needs to take each 2-handle separately and then see what happens when we collapse a 
1-handle it goes over using
the reference arc for that 1-handle. Let us give a few examples.

Start with the 2-handles in the $x-y$ plane, look at the pink 2-handle in that plane. When we use our system of 
reference arcs to collapse the 1-handles
it goes over we get that the pink 2-handle transforms to the following.

\centerline {\graphicspath{ {double_cover/Manifold_1011_double_cover/framing/x-y_plane/} } \includegraphics[width=11cm, height=8.5cm]{drawing_2}} 

As we are only interested in the framings of each 2-handle and not on how they link with the other 2-handles we have thrown away some 2-handles
so that the reader can easily see how the pink 2-handle will look like in dotted circle notation. It should be clear to the reader that
the pink 2-handle transforms to a knot with zero framing (remember we are using the blackboard framing).

The green 2-handle also transforms in a similar way, the following diagram shows how it looks like. The reader should notice that
the framing is again zero.

\centerline {\graphicspath{ {double_cover/Manifold_1011_double_cover/framing/x-y_plane/} } \includegraphics[width=11cm, height=8.5cm]{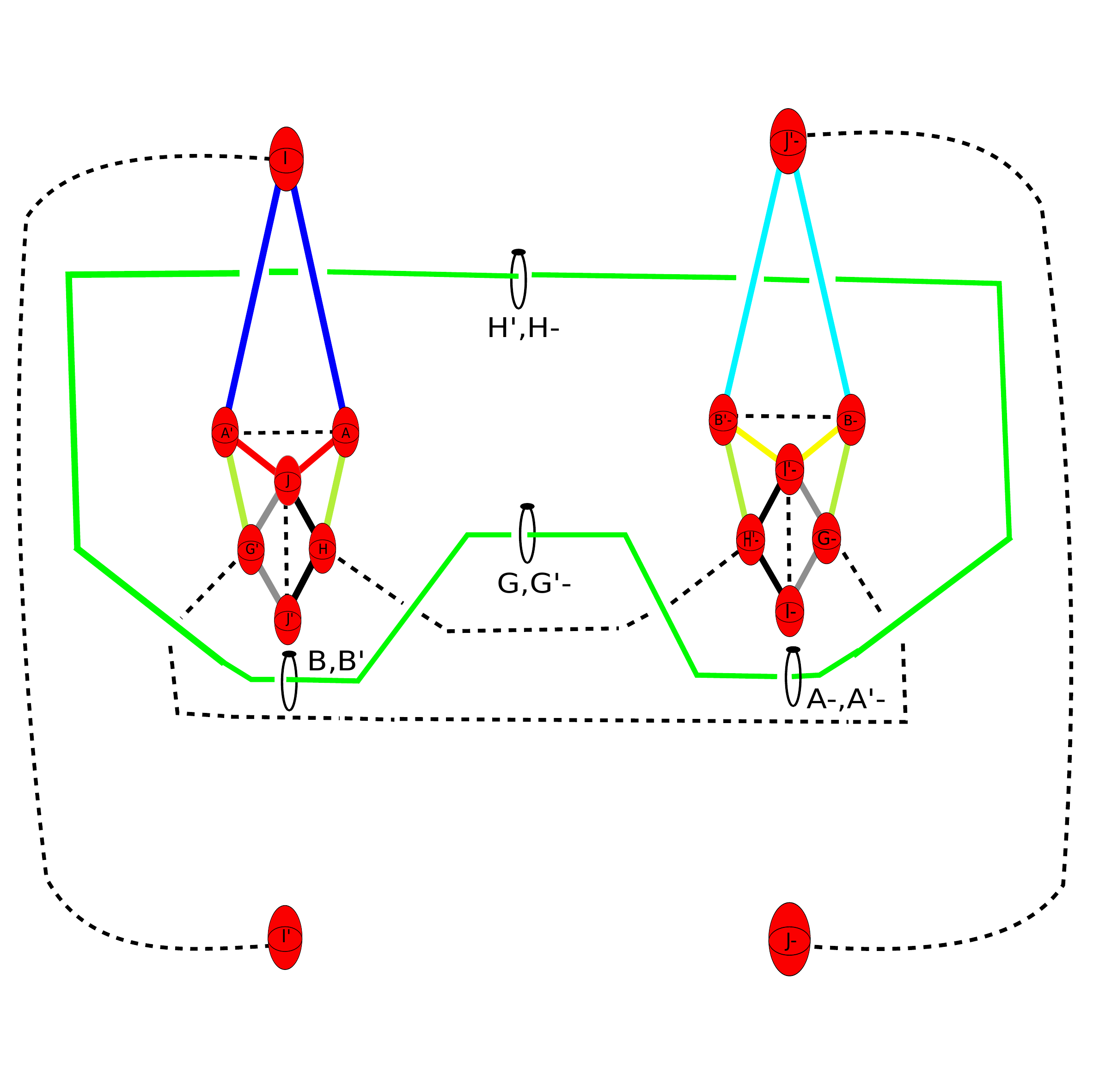}}

The following diagram shows how the red and blue 2-handles transform. It should be clear that the framing is zero for both of them.

\centerline {\graphicspath{ {double_cover/Manifold_1011_double_cover/framing/x-y_plane/} } \includegraphics[width=11cm, height=8cm]{drawing_1}}

One can continue in this fashion checking each 2-handle, in this case one finds that all the 2-handles will have framing zero. The point is that
since these 2-handles are all confined in a 2-plane one can easily find a system of reference arcs so that in dotted circle notation these 2-handles
will have even framing. A computation similar to the above works for the $x-z$ plane as well, we leave the details to the reader. We want to move
on to considering the diagram showing the twelve 2-handles that did not all lie in a single 2-plane.

The following diagram shows how the green 2-handle in that diagram transforms.

\centerline {\graphicspath{ {double_cover/Manifold_1011_double_cover/framing/x-y-z_plane/} } \includegraphics[width=12cm, height=10cm]{drawing_1}} 

It is easy to see that the framing is again zero. So far all the 2-handles we have shown have had framing zero simply because when they transform
they don't cross over themselves. It can happen that a 2-handle transforms to something that has self crossings. However, it is always the
case that the number of self crossings is even, this means that the blackboard framing will always be congruent to zero mod two. 

The following diagram shows how the black 2-handle transforms.

\centerline {\graphicspath{ {double_cover/Manifold_1011_double_cover/framing/x-y-z_plane/} } \includegraphics[width=12cm, height=10cm]{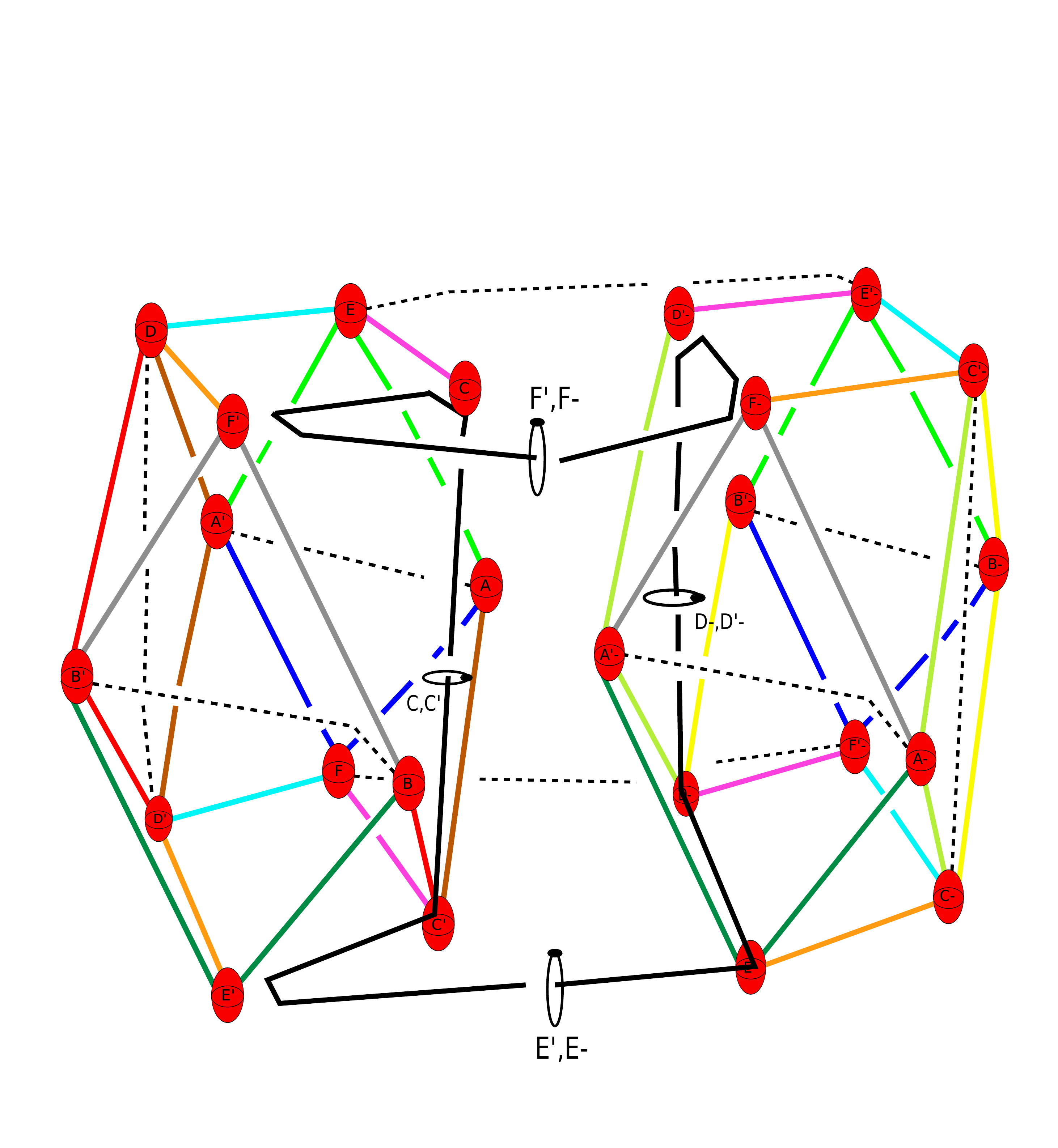}} 

In this case you can see two self crossings, and if we were to project this to the plane then we would have framing 
zero (as we count the crossings with sign).
Each of the other 2-handles in this diagram transform in a similar way. It is also the case that for the 2-handles 
in the diagram showing the $y-z$ plane that when we transform to dotted 
circle notation, using the reference arcs we have chosen above, we always get 2-handles with
an even number of self crossings. In other words they also all have an even framing. As the details are exactly 
analogous to what we have shown above we will not bother with the details.

In summary we have chosen a system of reference arcs that when used to transform the Kirby diagram of 
$\widetilde{M}$ into dotted circle notation we obtain 2-handles
all of which have an even framing. In particular by proposition \ref{spin} we see that the double cover of manifold 
1011 must be spin. 

We already mentioned that D. Ivan$\check{s}$i$\acute{c}$ in \cite{ivansic} proved that if you glue in a solid 3-
torus
to each cusp cross-section of $\widetilde{M}$ you get a 4-manifold that is homeomorphic to $S^4$.
The gluing is done by choosing a meridian, then gluing via the identity map. The parabolic translations
that corresponded to the $S^1$-fibre of each cusp cross-section of $M$ were $c$, $a$, $k$, $j$, and $e^{-1}g$.
Each of $c$, $a$, $k$, and $j$ are orientation preserving isometries, hence they lift to the same isometries
in the orientable double cover. This means we can take them as meridians in the orientable double cover.
On the other hand $e^{-1}g$ lifts to $(g^{-1}e)^{-1}(g^{-1}g)$ in the orientable double cover. On the level
of a handle decomposition the gluing in of a solid 3-torus corresponds to adding one 2-handle, two 3-handles, and
one 4-handle. Therefore on the level 
of the Kirby diagram of $\widetilde{M}$ the gluing in of solid 3-tori along the meridians
given by $c$, $a$, $k$, $j$, and $(g^{-1}e)^{-1}(g^{-1}g)$ corresponds to adding five 2-handles as follows.
We need to add a 2-handle between $C$ and $C'$ corresponding to gluing in a solid 3-torus via the identity along the 
meridian $c$. We then need to add a 2-handle from $A$ to $A'$ corresponding to gluing in a solid 3-torus via the 
identity along the meridian $a$, and similarly for $k$ and $j$. Finally we add a 2-handle that has two components
running from $E$ to $G$ and $g^{-1}G'$ to $g^{-1}E'$ respectively, corresponding to gluing in a solid 3-torus along
$(g^{-1}e)^{-1}(g^{-1}g)$. The Kirby diagram of $\widetilde{M}$ together with these five added 2-handles
constitutes a Kirby diagram for the filling of $\widetilde{M}$. Using proposition \ref{spin} and the system
of reference arcs we chose above, it is easy to see that each time we 
glue in one of these solid 3-tori we still obtain a spin 4-manifold. This observation will be important in the
next section when we construct certain finite covers of $\widetilde{M}$, and want to prove that their associated
fillings are all spin.

\section{Hyperbolic link complements in $\#_{2k}(S^2 \times S^2)$}\label{link_comp_1}

This section is devoted to constructing codimension two hyperbolic link complements in 4-manifolds that are 
homeomorphic to 
$\#_{2k}(S^2 \times S^2)$ for $k > 0$.
The candidate manifolds will be constructed by taking certain finite covers of manifold 1011 and then performing
a filling of these covers. The homeomorphism type of these fillings will be analysed using the fact that
the orientable double cover of manifold 1011 is spin (see section \ref{spin_str}) coupled with an argument showing
that the signature of these finite covers must vanish. As we will only be concerned with codimension two link
complements, from here on in we will simply use the term link complement.

Recall we are denoting manifold 1011 by $M$ and its orientable double cover by $\widetilde{M}$.
The manifold $M$ has five cusps with each cusp cross-section being given by the closed flat 3-manifold 
$\textbf{G}$. The parabolic translations corresponding to the $S^1$-fibre of each of these cusp cross-sections 
are given by $c$, $a$, $k$, $j$, $e^{-1}g$. When we construct a filling on the level of the fundamental group what 
we were doing is adding the relations
$c = a = k = j = e^{-1}g = 1$ (see end of section \ref{double_cover}). If instead we add a power of one of the 
translations then we can construct new fillings. For example
instead of taking $c$ we could take $c^n$ for some $n > 1$. This is still a translation and so we can try and 
analyse what happens if we fill
along this translation. On the level of the fundamental group what we are doing is adding the relations 
$c^n = a = k = j = e^{-1}g = 1$.
Using some computer software, such as \textbf{Magma}, one can obtain that the presentation for the fundamental group 
with these added relations
is $\langle c^n = e^2 = e^{-1}cec = 1 \rangle$. This is not quite a product as the $e$ and $c$ terms interact, 
however it is a semi direct product
of $\Z_n$, generated by $c$, and $\Z_2$, generated by $e$, where $e$ acts on $\Z_n$ by $e^{-1}ce$.  
We can use this fact to construct finite covers of $M$ as follows.

Let $\rho : \pi_1(M) \rightarrow  \Z_n \rtimes \Z_2$ be the quotient homomorphism given by the above discussion. We 
know explicit generators
of each parabolic subgroup associated to each cusp of $M$ (see section \ref{double_cover}). Using 
these explicit generators
one can show that the parabolic subgroups associated to the equivalence class of ideal vertices 
$(0, \pm 1, 0, 0)$, $(0, 0, \pm 1, 0)$, $(0, 0, 0, \pm 1)$, $(\pm 1/2, \pm 1/2, \pm 1/2, \pm 1/2)$ all map to the 
generator $e$ under $\rho$.
Furthermore, the parabolic subgroup corresponding to the ideal vertex $(\pm 1, 0, 0, 0)$ has one 
generator (other than $c$) that maps to $e$, since
$c$ maps to $c$ it follows that the image of the parabolic subgroup corresponding to $(\pm 1, 0, 0, 0)$ is the whole 
group $\Z_n \rtimes \Z_2$. 
This means that the image of the parabolic subgroup associated to $(\pm 1, 0, 0, 0)$ has index one, and the index of 
the image of all the other parabolic subgroups
is $n$. 
As the orientation reversing side pairing transformations of $M$ are $e$ and $g$, it follows that under the quotient 
map they map to $\Z_2$. This implies
that the orientation preserving generators are sent to $\Z_n$, which in turn implies that that orientable double 
cover of $M$ corresponds to the index
two subgroup given by $\rho^{-1}(\Z_n)$. We also know that the kernel of this quotient homomorphism is contained in 
$\rho^{-1}(\Z_n)$ and has index $2n$. It follows
that the cover corresponding to $ker(\rho)$ is also an $n$-fold cover of the orientable double cover of $M$. The 
associated deck transformation
group to $ker(\rho)$ over the orientable double cover of $M$ is $\Z_n$, which implies the $2n$ fold cover is a 
cyclic cover.

The above construction shows that the orientable double cover $\widetilde{M}$ has an $n$ sheeted cyclic cover for 
each $n$, or put another way $M$ has a $2n$ sheeted cover
for each $n$, where the case $n  = 1$ corresponds to the orientable double cover $\widetilde{M}$. The orientable 
double cover $\widetilde{M}$ has five cusps each having type $\textbf{A}$, which is the 3-torus.
The discussion above 
tells us how each boundary 3-torus corresponding
to a parabolic subgroup lifts to each cover. The boundary component corresponding to the ideal vertex class 
$(\pm 1, 0, 0, 0)$ lifts to one torus. This means
that the covering restricted to this boundary torus is the usual $n$ fold covering of the 3-torus onto itself. All 
other boundary 3-tori each lift
to $n$ copies of itself. In other words, the orientable double cover $\widetilde{M}$ has an $n$ sheeted cyclic 
cover which has $4n + 1$ cusps, each given by a 3-torus. We denote these covers by $\widetilde{M}_n$.

The question that arises at this point is can we identify any fillings of these cyclic covers? Observe that 
$\widetilde{M}$ has Euler
characteristic $2$ (as all the Ratcliffe-Tschantz manifolds have Euler characteristic $1$). Therefore each such 
cyclic cover has Euler characteristic $2n$.
In the case that $n$ is odd we claim that we can identify the topological type of the associated filling.
The argument proceeds as follows, from the previous section we know that $\widetilde{M}$ is spin hence it induces a 
spin structure on each of its boundary
3-tori (there are five in total). The three torus has exactly eight spin structures, seven of them are given by 
having one $S^1$-factor taking the
bounding spin structure induced from $D^2$, and each of these spin structures spin bounds a solid 3-torus.
The eighth spin structure corresponds to taking the Lie group spin structure on each $S^1$ factor. This
eighth spin structure does not spin bound a solid 3-torus, rather it spins bounds the complement of a singular 2-
torus 
fibre in $\CP^2\#_9\overline{\CP^2}$, viewed as an elliptic fibration over $\CP^1$ (see \cite{kirby_2} chap. V).
When we preform a filling of $\widetilde{M}$ we are gluing in a solid torus by making a choice of meridian. When we 
do this to each cusp cross-section the manifold we obtain is itself spin 
(see discussion at end of section \ref{spin_str}), meaning that the spin structure 
induced on each boundary 3-torus
must cross over the solid torus that we glued in. Hence the induced spin structure on each boundary component must 
be the one that spin bounds a solid torus, and
in fact the fibre we are filling in must have the spin structure induced from $D^2$. 
As four of the boundary 3-tori 
each lift to $n$ disjoint 3-tori in the cyclic
cover. It follows that the spin structure associated to these $n$ boundary 3-tori in the cyclic cover must have the 
same spin structure as their base.
Thus we need only worry about the lift of the boundary 3-torus corresponding to the ideal vertex class 
$(\pm 1, 0, 0, 0)$. In this case the covering
is given by the covering of $S^1 \times S^1 \times S^1$ onto itself induced by the $n$-fold cover of one $S^1$ 
factor onto itself. In the case that
$n$ is odd we have that the lifted spin structure on the associated boundary 3-torus in the $n$ sheeted cover must 
also be the one that spin bounds a
solid torus. This implies that for $n$ odd the associated filling of the $n$ sheeted cover is a closed smooth
simply connected spin 4-manifold. We have thus proved the following theorem.

\begin{thm}
The manifold $\widetilde{M}$ has a collection of $n$-sheeted covers, denoted by $\widetilde{M}_n$, such that
each $\widetilde{M}_n$ has $4n + 1$ cusps with each cusp cross-section being the 3-torus. Furthermore for 
$n$ odd a filling of $\widetilde{M}_n$ produces a closed smooth simply connected spin 4-manifold.
\end{thm}

We are not the first to construct the covers $\widetilde{M}_n$. In \cite{ivansic} p.18 Ivan$\check{s}$i$\acute{c}$
shows how to construct the manifolds $\widetilde{M}_n$ and states that they are complements in a 
closed smooth simply connected 4-manifold. However he does not identify the homeomorphism type of any of the
fillings of the $\widetilde{M}_n$.

We now want to move on to describe the homeomorphism type of each of the fillings of the
manifolds $\widetilde{M}_n$, in the case that $n > 1$ and odd. The starting point is to make use of the
following classification theorem of S. Donaldson and M. Freedman (see 
\cite{scorpan} p.244, see also \cite{donaldson} and \cite{freedman}).

\begin{thm}\label{simple_connec}
Every closed smooth simply connected 4-manifold is homeomorphic to either $S^4$ or $\#_m\CP^2\#_n\overline{\CP^2}$ 
or $\#_m\pm\mathcal{M}E_8\#_n(S^2\times S^2)$.
\end{thm}

The manifold $\mathcal{M}E_8$ denotes the non-smoothable simply connected 4-manifold with intersection form 
the $E_8$ lattice. The reader can find
a construction of this manifold and details about its intersection form in \cite{scorpan} p.86 and p.125. The 
manifold $-\mathcal{M}E_8$ denotes $\mathcal{M}E_8$ with opposite orientation.

The above classification theorem can be used to describe the homeomorphism type
of $\widetilde{M}_n$ for low values of $n$. For example when $n = 3$ the Euler characteristic 
of $\widetilde{M}_3$ is six, using the classification theorem  (theorem \ref{simple_connec}) we can deduce
that the filling of $\widetilde{M}_3$ (remember a filling does not change the Euler characteristic) 
must be homeomorphic to
one of $(S^2 \times S^2) \# (S^2 \times S^2)$, $\#_4\CP^2$, $\#_3\CP^2 \# \overline{\CP^2}$, $\#_2\CP^2 \#_2 
\overline{\CP^2}$, 
$\CP^2 \#_3 \overline{\CP^2}$ or $\#_4\overline{\CP^2}$. However only $(S^2 \times S^2) \# (S^2 \times S^2)$ is 
spin, it follows that $\widetilde{M}_3$ is homeomorphic to $(S^2 \times S^2) \# (S^2 \times S^2)$. 
We have thus proved the following proposition.

\begin{prop}
There exists a system $L$ consisting of thirteen linked tori embedded in a closed smooth 4-manifold $X$, such that 
$X$ is homeomorphic to $(S^2 \times S^2) \# (S^2 \times S^2)$, and
the complement $X - L$ admits a hyperbolic structure.
\end{prop}

When $n$ gets larger we get more homeomorphism types that admit spin structures. For example when
$n = 13$ we have that the filling of $\widetilde{M}_{13}$ has Euler characteristic $26$. In this case we have the
manifolds $K3 \# (S^2 \times S^2)$, where $K3$ denotes the complex algebraic surface with intersection form
$-2E_8 \oplus 3H$ ($H$ denotes the hyperbolic form), and the manifold $\#_{12}(S^2 \times S^2)$. The question
that arises is which one of  $K3 \# (S^2 \times S^2)$ or $\#_{12}(S^2 \times S^2)$ is the filling of
$\widetilde{M}_{13}$ homeomorphic to? The way to approach this problem is to start by understanding
the signature of the fillings of $\widetilde{M}_{n}$.  

The following theorem of D. Long and A. Reid allows us to compute the signature of any non-compact finite volume 
hyperbolic 4-manifold.

\begin{thm}\label{eta}
Let $M$ be a non-compact orientable finite volume hyperbolic 4-manifold. Then $\sigma(M) = \eta(\partial{M})$, where 
$\sigma$ denotes the signature and $\eta$ is the eta invariant.
\end{thm}

The proof of this theorem can be found in \cite{long} Theorem.2.1, p.173-174. We should mention that the above 
theorem is a special case
of a more general theorem proved by M.F. Atiyah, V.K. Patodi and I.M. Singer. Namely, they prove a general signature 
formula for a $4k$-dimensional
manifold with boundary see \cite{atiyah} Theorem.4.14, p.66. The proof of Long and Reid involves using the formula 
constructed by Atiyah, Patodi and Singer together with a theorem of Chern, which states that the Pontryagin classes 
of a hyperbolic manifold must vanish (\cite{chern}). 

The above theorem tells us that in order to work out the signature of the manifolds $\widetilde{M}_n$
we need to understand the eta invariant of the cusp cross-sections of $\widetilde{M}_n$. Each cusp
cross-section of $\widetilde{M}_n$ is a 3-torus, thus we need to compute the eta invariant of the
3-torus. The eta invariant of the flat 3-torus is known to be zero, the computation of this can be
found in \cite{szcz} Example 1, p.128. This gives the following proposition.

\begin{prop}
The hyperbolic manifolds $\widetilde{M}_n$ have vanishing signature. 
\end{prop}

When we carry out a filling we are gluing in a solid 3-torus by filling in one of the
$S^1$-fibres with a copy of $D^2$ (i.e. we are doing a filling with a fixed choice of meridian). A solid 3-torus 
also has vanishing signature, therefore when we perform a filling we are gluing two 4-manifolds with signature zero
along their common boundary (via the identity map).
There is a beautiful theorem of S. Novikov that allows one to compute the signature of a gluing provided one knows 
the signature of the pieces that are being glued together (see \cite{kirby_2} Thm.5.3, p.27).

\begin{thm}\label{novi}
Given two oriented $4n$-dimensional manifolds $M$ and $N$ such that $\partial{M} = \partial{N}$. 
Then $\sigma(M \cup_{\partial} N) = \sigma(M) + \sigma(N)$, 
where $M \cup_{\partial} N$ denotes $M$ glued to $N$ along the common boundary.  
\end{thm}

This theorem implies that the filling of any of the $n$ sheeted covers of $\widetilde{M}_n$ must have signature 
zero. Furthermore for $n$ odd we know that
such a filling must be spin. This observation coupled with the classification theorem of S. Donaldson and M. 
Freedman allows us to completely describe the homeomorphism type of the fillings of 
$\widetilde{M}_n$ for $n$ odd.

\begin{prop}\label{mainthm_3}
For $n > 1$ odd the $n$ sheeted cover $\widetilde{M}_n$ is a complement of $4n+1$ tori in a manifold $X$ that is 
homeomorphic to $\#_{n-1}(S^2 \times S^2)$.
\end{prop}

These examples give the following theorem.

\begin{thm}
For $k \geq 1$ there exists a collection of $8k + 5$ tori embedded in a closed smooth 4-manifold $X$, such that
$X$ is homeomorphic to $\#_{2k}(S^2 \times S^2)$, and $X-L$ admits a hyperbolic structure.
\end{thm}

We remark that in the case that $n$ is even the bounding spin structure on the base $S^1$ can lift to the 
non-bounding (the Lie group spin structure)
spin structure on the total space $S^1$. It is this that obstructs us from being able to conclude that the $n$ 
sheeted cover of $\widetilde{M}$ for $n$ even is spin.

\section{Hyperbolic link complements in closed smooth simply connected 4-manifolds}\label{link_comp_2}

In the previous section we witnessed the construction of hyperbolic link complements in closed 4-manifolds that were 
homeomorphic to $\#_{2k}(S^2 \times S^2)$, for $k > 0$. The identification of the homeomorphism type proceeded via 
three main steps, which can be broken down as follows. The first step involved using the
classification theorem of S. Donaldson and M. Freedman (see theorem \ref{simple_connec}) to give a finite list of 
homeomorphism types, for each 
Euler characteristic, from which we could choose from. The second step involved refining this list to have fewer 
possible 
choices by proving that the fillings we were dealing with
had to admit a spin structure. The third and final step involved a further refinement by the understanding that the 
fillings also had to have vanishing signature. The need for the second step was primarily to do with the fact of 
understanding whether the manifolds we were dealing with
could have the topological type of $\#_r\CP^2  \#_s\overline{\CP^2}$. As the manifolds we were dealing with were all 
spin, and all the manifolds
$\#_r\CP^2  \#_s\overline{\CP^2}$ are not, we could cross these off our list. This brings forth the question of 
whether a manifold that is homeomorphic to $\#_r\CP^2  \#_s\overline{\CP^2}$ could have a link complement that is 
hyperbolic? More generally we could ask which closed smooth simply connected 4-manifolds can admit a 
hyperbolic link complement? In general this question is very hard to answer, the main reason being that
our understanding of four dimensional hyperbolic link complements is still in its infancy. In this brief
section we give necessary conditions on the homeomorphism type for a closed smooth simply connected 4-manifold
to contain a link complement that admits a hyperbolic structure.

Recall from section \ref{link_comp_1} that a theorem of D. Long and A. Reid (theorem \ref{eta}) tells us
that the signature of a finite volume non-compact hyperbolic 4-manifold is given by the eta
invariant of the cusp cross-sections.
The importance of this theorem is that it pushes the computation of the signature to the computation of 
the eta invariant of the cusp cross-sections of the
manifold. The cusp cross-sections of a non-compact orientable finite volume hyperbolic 4-manifold fall into six 
classes, coming from the fact that
there are six isometry classes of orientable closed flat 3-manifolds. We denoted these six classes by $\textbf{A}$, 
$\textbf{B}$,
$\textbf{C}$, $\textbf{D}$, $\textbf{E}$, $\textbf{F}$ 
(or in Wolf's notation by $\mathcal{G}_1$, $\mathcal{G}_2$, $\mathcal{G}_3$, $\mathcal{G}_4$, $\mathcal{G}_5$, 
$\mathcal{G}_6$). Therefore in order
to understand the signature of a non-compact orientable finite volume hyperbolic 4-manifold one needs to understand 
what the eta invariant of
the above six classes of flat 3-manifolds are.    

The computation of the eta invariant for these six classes of flat 3-manifolds can be found in \cite{szcz} Example 
1, p.128. The following proposition gives the values of the eta invariant for these six classes.

\begin{prop}\label{eta_flat}

\item $\eta(\textbf{A}) = 0$

\item $\eta(\textbf{B}) = 0$

\item $\eta(\textbf{C}) = \frac{-2}{3}$

\item $\eta(\textbf{D}) = -1$

\item $\eta(\textbf{E}) = \frac{-4}{3}$

\item $\eta(\textbf{F}) = 0$

\end{prop}

From the
classification theorem (see \cite{wolf} Theorem.3.5.5, p.117) it is known that only $\textbf{A}$ and $\textbf{B}$ 
are $S^1$-fibre bundles
over a compact surface, with $\textbf{A}$ being a 3-torus fibering over a 2-torus and $\textbf{B}$ fibering over a 
Klein bottle. As we are only
interested in link complements (remember link complement for us means codimension two link complement) we need only 
worry about the eta invariant of $\textbf{A}$ and $\textbf{B}$. Using 
the above proposition
we obtain the following corollary.

\begin{cor}
Let $M$ be an orientable non-compact finite volume hyperbolic 4-manifold with cusp cross-sections of type 
$\textbf{A}$ or $\textbf{B}$. Then $\sigma(M) = 0$, where $\sigma$ is the signature invariant.
\end{cor}

From now on we assume $M$ is an orientable non-compact finite volume hyperbolic 4-manifold with cusp cross-sections 
of type $\textbf{A}$ or $\textbf{B}$.
When we construct a filling of $M$ we are gluing in disk bundles associated to the cusp cross-sections of type 
$\textbf{A}$ or $\textbf{B}$.
These disk bundles also have vanishing signature, hence the filling of $M$ is made up from a collection of 
4-manifolds each having vanishing signature.
Appealing to Novikov's theorem (see theorem \ref{novi})
we can conclude that the signature of such a filled in manifold must be zero. In other words if
a closed smooth 4-manifold admits a hyperbolic link complement it must have have vanishing signature.
This simple observation on the signature gives necessary conditions on the homeomorphism type of those closed
smooth simply connected 4-manifolds that admit a hyperbolic link complement.

\begin{thm}\label{mainthm_5}
Let $M$ be a closed smooth simply connected 4-manifold that has a link complement that admits  
hyperbolic structure. Then the homeomorphism
type of $M$ falls into one of the following three categories:

\begin{itemize}
\item $S^4$

\item $\#_k(S^2 \times S^2)$, $k > 0$.

\item $\#_k\CP^2 \#_k \overline{\CP^2}$, $k > 0$.
\end{itemize}
\end{thm}

Using this theorem we obtain the following corollary.

\begin{cor}
Let $M$ be a closed irreducible smooth simply connected 4-manifold that contains a link complement that admits
a hyperbolic structure. The $M$ is homeomorphic to $S^4$ or $S^2 \times S^2$.
\end{cor} 

In the irreducible case examples of link complements that admit a hyperbolic structure have already
been constructed. For 4-manifolds that are homeomorphic to $S^4$ this was first done by
D. Ivan$\check{s}$i$\acute{c}$ in \cite{ivansic}, then later on by 
D. Ivan$\check{s}$i$\acute{c}$, J. Ratcliffe, and S. Tschantz in \cite{tschantz}. The case of an example in a 4-
manifold
with homeomorphism type $S^2 \times S^2$ was done by the author in \cite{sarat_3} (we actually showed that
our example was diffeomorphic to a standard $S^2 \times S^2$).
For the reducible case, we saw in the previous section
the construction of examples in $\#_{2k}(S^2 \times S^2)$. In \cite{sarat_3} we constructed
a hyperbolic link complement in a manifold that is diffeomorphic to $S^2 \times S^2$, we believe one can use this
manifold and a similar construction to what we did for the case of $\#_{2k}(S^2 \times S^2)$ to produce
examples with homeomorphism type $\#_{2k+1}(S^2 \times S^2)$. This leaves the question of
constructing hyperbolic link complements in manifolds with homeomorphism type $\#_k\CP^2 \#_k \overline{\CP^2}$.

\textbf{Question:} Does there exist a hyperbolic link complement in $\#_k\CP^2 \#_k \overline{\CP^2}$ for
some $k > 0$?

At present we do not know of any candidate hyperbolic 4-manifold that could lead to a construction
of a hyperbolic link complement in $\#_k\CP^2 \#_k \overline{\CP^2}$. 
The four dimensional hyperbolic link 
complements that have been constructed to this day have all made use of the hyperbolic 4-manifolds constructed
by Ratcliffe and Tschantz in \cite{ratcliffe}. Their census contains 1171 distinct isometry classes of 
finite volume non-compact hyperbolic 4-manifolds and it may be possible that, just as before, one of these
could lead to the construction of a hyperbolic link complement in $\#_k\CP^2 \#_k \overline{\CP^2}$.


\begin{thebibliography}{10}

\bibitem{atiyah}
Atiyah, M.F. \emph{Spectral asymmetry and Riemannian 
Geometry. I}, Math. Proc. Camb. Phil. Soc. (1975), 77, 43.


\bibitem{chern}
Chern, S. \emph{On curvature and characteristic classes of a Riemannian 
manifold}, Abhandlungen aus dem Mathematischen Seminar der Universität Hamburg (1955), 20.


\bibitem{donaldson}
Donaldson, S.K. \emph{An application of gauge theory to four dimensional  
topology}, J. Differential Geometry, volume 18, number 2 (1983), 279-315.


\bibitem{freedman}
Freedman, M.H. \emph{The topology of four-dimensional 
manifolds}, J. Differential Geometry, volume 17, number 3 (1982), 357-453.

\bibitem{gompf}
Gompf, R.E. and Stipsicz, A.I \emph{4-manifolds and Kirby Calculus
}, Graduate Studies in Mathematics, Providence, Rhode Island, 1999.


\bibitem{hantzsche}
Hantzsche, W. and Wendt, H. \emph{Dreidimensionale euklidische Raumformen
}, Math. Ann. 110 (1935), 593-611.

\bibitem{ivansic}
Ivan$\check{s}$i$\acute{c}$, D. \emph{Hyperbolic structure on a 
complement of tori in the 4-sphere}, Adv. Geom. 4, no. 1 (2004), 119–139.

\bibitem{tschantz}
Ivan$\check{s}$i$\acute{c}$, D. and Ratcliffe, J.G. and Tschantz, S.T. \emph{Complements of tori and Klein bottles 
in the 4-sphere that have hyperbolic structure}, Algebr. Geom. Topol. 5 (2005), 999–1026.


\bibitem{kirby_2}
Kirby, R. \emph{Topology of
4-manifolds}, Springer-Verlag, Berlin Heidelberg, 1989.

\bibitem{long}
Long, D. and Reid, A. \emph{On the geometric boundaries of hyperbolic
4-manifolds}, Geometry \& Topology, Vol. 4 (2000) no. 5, p.171-178.


\bibitem{ratcliffe}
Ratcliffe, J.G. and Tschantz, S.T. \emph{The Volume Spectrum of
 Hyperbolic 4-manifolds}, Experiment. Math.
Volume 9, Issue 1 (2000), 101-125.

\bibitem{sarat}
Saratchandran, H. \emph{Kirby diagrams and the Ratcliffe-Tschantz 
hyperbolic 4-manifolds}, ArXiv:math:GT/1503.06722.

\bibitem{sarat_2}
Saratchandran, H. \emph{A four dimensional hyperbolic link complement 
in a standard $S^4$}, ArXiv:math:GT/1503.07778.

\bibitem{sarat_3}
Saratchandran, H. \emph{A four dimensional hyperbolic link complement 
in a standard $S^2 \times S^2$}, ArXiv:math:GT/1504.01366.

\bibitem{scorpan}
Scorpan, A. \emph{The Wild World of 
4-manifolds}, American Mathemtical Society, Providence, Rhode island, 2005. 


\bibitem{szcz}
Szczepa$\acute{n}$ski, A. \emph{Eta invariants for flat 
manifolds}, Annals of Global Analysis and Geometry, February 2012, Volume 41, Issue 2, pp 125-138.



\bibitem{wolf}
Wolf, J.A. \emph{Spaces of constant 
curvature}, McGraw-Hill, United States of America, 1967.

\end{thebibliography}
\end{document}